\input boxedeps.tex % unix
\SetepsfEPSFSpecial % unix
\HideDisplacementBoxes
\def\figin#1#2{\medbreak
$$
 {\BoxedEPSF{#1 scaled
#2}%
}%
$$
\medbreak\noindent}
\documentstyle{amltd}
\begin{document}

\annalsline{151}{2000}
\received{January 4, 1999}
\startingpage{327}
\def\bye{\end{document}}
 \font\tenrm=cmr10
\input amssym.def
\input amssym.tex

\newcommand{\sumg}{\!\!\sum_{\ell_\gamma(X) < \varepsilon}}
\newcommand{\sumd}{\!\!\sum_{\ell_\delta(X) < 2\varepsilon}}
\newcommand{\sumdg}{\!\!\sum_{\varepsilon \le \ell_\delta(X) < 2\varepsilon}}
%for Bbb in amstex
\catcode`\@=11
\font\twelvemsb=msbm10 scaled 1100
\font\tenmsb=msbm10
%\font\ninemsb=msbm7 scaled 1100%msbm9
\font\ninemsb=msbm10 scaled 800
\newfam\msbfam
\textfont\msbfam=\twelvemsb  \scriptfont\msbfam=\ninemsb
  \scriptscriptfont\msbfam=\ninemsb
\def\msb@{\hexnumber@\msbfam}
\def\Bbb{\relax\ifmmode\let\next\Bbb@\else
 \def\next{\errmessage{Use \string\Bbb\space only in math
mode}}\fi\next}
\def\Bbb@#1{{\Bbb@@{#1}}}
\def\Bbb@@#1{\fam\msbfam#1}
\catcode`\@=12

 \catcode`\@=11
\font\twelveeuf=eufm10 scaled 1100
\font\teneuf=eufm10
\font\nineeuf=eufm7 scaled 1100%eufm9
\newfam\euffam
\textfont\euffam=\twelveeuf  \scriptfont\euffam=\teneuf
  \scriptscriptfont\euffam=\nineeuf
\def\euf@{\hexnumber@\euffam}
\def\frak{\relax\ifmmode\let\next\frak@\else
 \def\next{\errmessage{Use \string\frak\space only in math
mode}}\fi\next}
\def\frak@#1{{\frak@@{#1}}}
\def\frak@@#1{\fam\euffam#1}
\catcode`\@=12

%--------------- Author macros ---------------

\input geompsfi.sty

 %-----------------------------------------------------------
%
%	General Math Macros
%
%-----------------------------------------------------------
%	\usepackage{amsmath, amssymb, latexsym}
%-----------------------------------------------------------
%
%  Fonts synonyms
%
\newcommand{\gothic}{\frak}
\newcommand{\bb}{\Bbb}
%	Blackboard symbols
%
\newcommand{\aff}{{\bb A}}
\newcommand{\ball}{{\bb B}}
\newcommand{\cx}{{\bb C}}
\newcommand{\eucl}{{\bb E}}
\newcommand{\field}{{\bb F}}
\newcommand{\half}{{\bb H}}
\newcommand{\integers}{{\bb Z}}
\newcommand{\lowerhalf}{{\bb L}}
\newcommand{\natls}{{\bb N}}
\newcommand{\ratls}{{\bb Q}}
\newcommand{\reals}{{\bb R}}
\newcommand{\proj}{{\bb P}}
%
%  Figures with captions and labels
% \makefig{Title.}{fig:label}{contents}
%
%\newcommand{\makefig}[3]{
%	\begin{figure}[htbp]
%        \refstepcounter{figure}
%	\label{#2}
 %       \begin{center}
	%	~#3~\\
	%	\medskip
%                {\sf Figure \thefigure.  #1}
 %       \end{center}
	%\end{figure}
%}
%
% For frames in figures
%
\thicklines
%
% Recalling theorems already stated.
% It is imperitive to issue a recallend after
% *each* recallthm.
%
\newcommand{\recallthm}[1]{
  \renewcommand{\thetheorem}{\ref{#1}}
}
 
\newcommand{\recallend}{
 \addtocounter{theorem}{-1}
 \renewcommand{\thetheorem}{\thechapter.\arabic{theorem}}
}

%
%  Tables with captions and labels
%
\newcommand{\maketab}[3]{
	\begin{figure}[htbp]
        \refstepcounter{figure}
	\label{#2}
        \begin{center}
		#3~\\
		\bigskip
                {\sf Table \thefigure.  #1}
        \end{center}
	\end{figure}
}
%
%  Additional table features:
%  Add some room; use "\roll" for multirow entries
%
\renewcommand{\arraystretch}{1.25}
\newcommand{\roll}{\vspace{-.4em}\\}

%         Lengths
\newcommand{\widemargins}{
\setlength{\textwidth}{6.0in}
\setlength{\oddsidemargin}{0.25in}
\setlength{\evensidemargin}{0.25in}
}

\newcommand{\highmargins}{
\setlength{\headheight}{0.0in}
\setlength{\headsep}{0.0in}
\setlength{\topmargin}{0.0in}
\setlength{\footskip}{0.0in}
\setlength{\textheight}{9.5in}
}

%	Labelled negatively indented paragraphs.
%	Usage: \begin{negpar}
%	       \item
%	Paragraph whose first line
%		only, is negatively indented.
%		\end{negpar}
\newenvironment{negpar}{
        \begin{list}{}{ \setlength{\leftmargin}{.5in}
                        \setlength{\itemindent}{-\leftmargin}
                        }
                }{
        \end{list}
                } 
%
%	Emphasized paragraph label.
%
 
\newcommand{\mat}[1]{\left( \begin{smallmatrix} #1 \end{smallmatrix}\right) }
\newcommand{\Mat}[1]{\left(\begin{matrix} #1 \end{matrix}\right)}

% 
% Wide tildes standard
%
\renewcommand{\tilde}{\widetilde}

%
% Under-arrow projective/direct limits
%  (instead of text "projlim", "dirlim")
%
%\renewcommand{\projlim}{\underset{\longleftarrow}{\lim\;}}
%\newcommand{\dirlim}{\underset{\longrightarrow}{\lim\;}}

%      Math Synonyms
\newcommand{\actson}{\leadsto}
\newcommand{\arrow}{\rightarrow}
\newcommand{\Arrow}{\longrightarrow}
% Asymptotic to
\newcommand{\asyto}{\sim}
% Bounded Above and Below
\newcommand{\bab}{\asymp}
\newcommand{\bdry}{\partial}
\newcommand{\bijects}{\leftrightarrow}
\newcommand{\birat}{\dashrightarrow}
\newcommand{\brackets}[1]{\langle #1 \rangle}
\newcommand{\Brackets}[1]{\left\langle #1 \right\rangle}
\newcommand{\ceiling}[1]{\lceil #1 \rceil}
\newcommand{\closure}{\overline}
\newcommand{\compos}{\circ}
\newcommand{\compactsubset}{\Subset}
\newcommand{\congruent}{\equiv}
\newcommand{\conjugate}{\overline}
\newcommand{\contains}{\supset}
\newcommand{\degree}{^{\circ}}
\newcommand{\degrees}{^{\circ}}
\newcommand{\del}{\partial}
\newcommand{\dirminus}{\ominus}
\newcommand{\dirsum}{\oplus}
\newcommand{\Dirsum}{\bigoplus}
\newcommand{\disjunion}{\sqcup}
\newcommand{\bigdisjunion}{\amalg}
\newcommand{\equi}{\sim}
\newcommand{\equii}{\approx}
\newcommand{\grad}{\nabla}
%	Already present in amslatex.
 \newcommand{\impliedby}{\Longleftarrow}
 \newcommand{\implies}{\Longrightarrow}
\newcommand{\includesin}{\hookrightarrow}
\newcommand{\invlim}[1]{\lim_{\stackrel{\textstyle \longleftarrow}{#1}}}
\newcommand{\isom}{\cong}
\newcommand{\join}{\vee}
%	Already present in amslatex.
%\newcommand{\iff}{\Longleftrightarrow}
\newcommand{\lap}{\Delta}
\newcommand{\mapsonto}{\twoheadrightarrow}
\newcommand{\mem}{\in}
\newcommand{\notarrow}{\nrightarrow}
\newcommand{\notmem}{\not\in}
\newcommand{\Not}{\sim}
\newcommand{\nullset}{\emptyset}
\newcommand{\nsubgp}{\rhd}
\newcommand{\plusorminus}{\pm}
\newcommand{\symdif}{\triangle}
\newcommand{\superset}{\supset}
\newcommand{\tensor}{\otimes}
\newcommand{\wave}{\square}

%      Math Spacing Improvements
\newcommand{\BAB}{\;\;\bab\;\;}
\newcommand{\LE}{\;\;\le\;\;}
\newcommand{\GE}{\;\;\ge\;\;}
\newcommand{\LT}{\;\;<\;\;}
\newcommand{\GT}{\;\;>\;\;}
\newcommand{\EQ}{\;\;=\;\;}
\newcommand{\dx}{\;dx}
\newcommand{\st}{\; : \;}         %Such that
\newcommand{\prodprime}{\sideset{}{'}\prod}
\newcommand{\sumprime}{\sideset{}{'}\sum}

%      Math Overlines
\newcommand{\Abar}{{\overline{A}}}
\newcommand{\Bbar}{{\overline{B}}}
\newcommand{\Cbar}{{\overline{C}}}
\newcommand{\Dbar}{{\overline{D}}}
\newcommand{\Kbar}{{\overline{K}}}
\newcommand{\Mbar}{{\overline{M}}}
\newcommand{\Nbar}{{\overline{N}}}
\newcommand{\Xbar}{{\overline{X}}}
\newcommand{\Sbar}{{\overline{S}}}
\newcommand{\Tbar}{{\overline{T}}}
\newcommand{\Ubar}{{\overline{U}}}
\newcommand{\Ybar}{{\overline{Y}}}
\newcommand{\Zbar}{{\overline{Z}}}
\newcommand{\abar}{{\overline{a}}}
\newcommand{\bbar}{{\overline{b}}}
\newcommand{\cbar}{{\overline{c}}}
\newcommand{\dbar}{{\overline{d}}}
\newcommand{\ebar}{{\overline{e}}}
\newcommand{\fbar}{{\overline{f}}}
\newcommand{\gbar}{{\overline{g}}}
\newcommand{\kbar}{{\overline{k}}}
\newcommand{\rbar}{{\overline{r}}}
\newcommand{\sbar}{{\overline{s}}}
\newcommand{\tbar}{{\overline{t}}}
\newcommand{\ubar}{{\overline{u}}}
\newcommand{\vbar}{{\overline{v}}}
\newcommand{\wbar}{{\overline{w}}}
\newcommand{\xbar}{{\overline{x}}}
\newcommand{\ybar}{{\overline{y}}}
\newcommand{\zbar}{{\overline{z}}}

\newcommand{\alphabar}{\overline{\alpha}}
\newcommand{\dwbar}{{d \wbar}}
\newcommand{\dzbar}{{d \zbar}}
\newcommand{\delbar}{{\overline{\del}}}
\newcommand{\Deltabar}{{\overline{\Delta}}}
\newcommand{\Gammabar}{{\overline{\Gamma}}}
\newcommand{\mubar}{{\overline{\mu}}}
\newcommand{\Omegabar}{{\overline{\Omega}}}
\newcommand{\omegabar}{{\overline{\omega}}}
\newcommand{\phibar}{{\overline{\phi}}}
\newcommand{\Phibar}{{\overline{\Phi}}}
\newcommand{\psibar}{{\overline{\psi}}}

%	Hats
\newcommand{\ahat}{{\widehat{a}}}
\newcommand{\chat}{{\widehat{\cx}}}
\newcommand{\fhat}{{\widehat{f}}}
\newcommand{\ghat}{{\widehat{g}}}
\newcommand{\hhat}{{\widehat{h}}}
\newcommand{\rhat}{{\widehat{\reals}}}
\newcommand{\Ahat}{{\widehat{A}}}
\newcommand{\Ghat}{{\widehat{G}}}
\newcommand{\Hhat}{{\widehat{H}}}
\newcommand{\Khat}{{\widehat{K}}}
\newcommand{\Nhat}{{\widehat{N}}}
\newcommand{\Phat}{{\widehat{P}}}
\newcommand{\Qhat}{{\widehat{Q}}}
\newcommand{\rhohat}{{\widehat{\rho}}}
\newcommand{\psihat}{{\widehat{\psi}}}
\newcommand{\uhat}{{\widehat{u}}}
\newcommand{\vhat}{{\widehat{v}}}
\newcommand{\xhat}{{\widehat{x}}}
\newcommand{\yhat}{{\widehat{y}}}
\newcommand{\zedhat}{{\widehat{\zed}}}

%	Math Words
\newcommand{\ad}{{\rm ad}}
\newcommand{\Ad}{{\rm Ad}}
\newcommand{\area}{{\rm area}}
\newcommand{\Aut}{{\rm Aut}}
\newcommand{\avg}{{\rm avg}}
\newcommand{\boxdim}{{\rm box-dim}}
\newcommand{\Coker}{{\rm Coker}}
\newcommand{\core}{{\rm core}}
\newcommand{\Cov}{{\rm Cov}}
\newcommand{\CS}{{\rm CS}}
\newcommand{\curl}{{\rm curl}}
\newcommand{\diag}{{\rm diag}}
\newcommand{\diam}{{\rm diam}}
\newcommand{\Diff}{{\rm Diff}}
\newcommand{\dist}{{\rm dist}}
\renewcommand{\div}{{\rm div}}
\newcommand{\Div}{{\rm Div}}
\newcommand{\E}{{\rm E}}
\newcommand{\End}{{\rm End}}
\newcommand{\Gal}{{\rm Gal}}
\newcommand{\gl}{{\rm gl}}
\newcommand{\GL}{{\rm GL}}
\newcommand{\gr}{{\rm gr}}
\newcommand{\Gr}{{\rm Gr}}
\newcommand{\Hdim}{{\rm H.dim}}
\newcommand{\Hom}{{\rm Hom}}
\newcommand{\hull}{{\rm hull}}
\newcommand{\hyp}{{{\rm hyp}}}
\newcommand{\Inner}{{\rm Inner}}
\newcommand{\id}{{\rm id}}
\newcommand{\ind}{{\rm ind}}
\newcommand{\interior}{{\rm int}}
\renewcommand{\Im}{{\rm Im}}
\newcommand{\Ind}{{\rm Ind}}
\newcommand{\Isom}{{\rm Isom}}
\newcommand{\Jac}{{\rm Jac}}
\newcommand{\Ker}{{\rm Ker}}
\newcommand{\lcm}{{\rm lcm}}
\newcommand{\length}{{\rm length}}
\newcommand{\lk}{{\rm lk}}
\newcommand{\logplus}{{\rm \stackrel{+}{\log}}}
\newcommand{\ml}{{\cal M}{\cal L}}
% Renew Required For amslatex
\def\mod{{\rm mod}}
\newcommand{\Mod}{{\rm Mod}}
\newcommand{\mult}{{\rm mult}}
\newcommand{\ord}{{\rm ord}}
\newcommand{\Out}{{\rm Out}}
\newcommand{\Pic}{{\rm Pic}}
\newcommand{\Poly}{{\rm Poly}}
\newcommand{\Proj}{{\rm Proj}}
\newcommand{\PSL}{{\rm PSL}}
\newcommand{\rad}{{{\rm rad}}}
\newcommand{\rank}{{\rm rank}}
\newcommand{\Rat}{{\rm Rat}}
\renewcommand{\Re}{{\rm Re}}
\newcommand{\Res}{{\rm Res}}
\newcommand{\Rinf}{\reals_\infty}
\newcommand{\Roots}{{\rm Roots}}
\newcommand{\sign}{{\rm sign}}
\newcommand{\Sinf}{S_\infty}
\newcommand{\Sinfn}{{S_\infty^{n-1}}}
\renewcommand{\sl}{{\rm sl}}
\newcommand{\SL}{{\rm SL}}
\newcommand{\so}{{\rm so}}
\newcommand{\SO}{{\rm SO}}
\newcommand{\Sp}{{\rm Sp}}
\newcommand{\Spec}{{\rm Spec}}
\newcommand{\Stab}{{\rm Stab}}
\newcommand{\supp}{{\rm supp}}
\newcommand{\Teich}{{\rm Teich}}
\newcommand{\Tr}{{\rm Tr}}
\newcommand{\tr}{{\rm tr}}
\newcommand{\Var}{{\rm Var}}
\newcommand{\vol}{{\rm vol}}
\newcommand{\zed}{\integers}

% Math words from amssymb package
\newcommand{\semidir}{\ltimes}

%   \cX  = draw X caligraphically.
\newcommand{\cA}{{\cal A}}
\newcommand{\cB}{{\cal B}}
\newcommand{\cC}{{\cal C}}
\newcommand{\cD}{{\cal D}}
\newcommand{\cE}{{\cal E}}
\newcommand{\cF}{{\cal F}}
\newcommand{\cG}{{\cal G}}
\newcommand{\cH}{{\cal H}}
\newcommand{\cI}{{\cal I}}
\newcommand{\cJ}{{\cal J}}
\newcommand{\cK}{{\cal K}}
\newcommand{\cL}{{\cal L}}
\newcommand{\cM}{{\cal M}}
\newcommand{\cN}{{\cal N}}
\newcommand{\cO}{{\cal O}}
\newcommand{\cP}{{\cal P}}
\newcommand{\cQ}{{\cal Q}}
\newcommand{\cR}{{\cal R}}
\newcommand{\cS}{{\cal S}}
\newcommand{\cT}{{\cal T}}
\newcommand{\cU}{{\cal U}}
\newcommand{\cV}{{\cal V}}
\newcommand{\cW}{{\cal W}}
\newcommand{\cX}{{\cal X}}
\newcommand{\cY}{{\cal Y}}
\newcommand{\cZ}{{\cal Z}}

%   Math names
\newcommand{\Cech}{\v{C}ech~}
\newcommand{\Mane}{Ma{\~n}\'e~}
\newcommand{\Swiatek}{\'{S}wi\c{a}tek~}

\newcommand{\T}{{\rm T}}
\newcommand{\Log}{{\rm Log}}
\newcommand{\WP}{{{\rm WP}}}

\newcommand{\od}{{1/d}}
\newcommand{\god}{g_\od}
\newcommand{\wod}{\omega_\od}

\newcommand{\ol}{{1/\ell}}
\newcommand{\gol}{g_\ol}
\newcommand{\wol}{\omega_\ol}
\newcommand{\tol}{\theta_\ol}

\newcommand{\sigf}{\sigma_F}
\newcommand{\sigqf}{\sigma_{QF}}
\newcommand{\sigqfbar}{\overline{\sigma}_{QF}}

%-------------- Author entries --------------------

\title{The moduli space of Riemann surfaces\\ is K\"ahler hyperbolic} %Article title
\shorttitle{Moduli space of Riemann surfaces}  
  \acknowledgements{Research partially supported by the NSF.\hfill\break
\phantom{1}\hskip.19in
        1991 Mathematics Subject Classification:
        Primary 32G15, Secondary 20H10, 30F60, 32C17.}
 \author{Curtis T. McMullen}
 \institutions{Harvard University, Cambridge, MA}

\bigbreak
\centerline{\bf Contents}
\medbreak
\noindent 1\quad  Introduction

\noindent 2\quad  Teichm\"uller space

\noindent 3\quad  $1/\ell$ is almost pluriharmonic

\noindent 4\quad  Thick-thin decomposition of quadratic differentials

\noindent 5\quad  The $1/\ell$ metric

\noindent 6\quad  Quasifuchsian reciprocity

\noindent 7\quad  The Weil-Petersson form is $d$(bounded)

\noindent 8\quad  Volume and curvature of moduli space

\noindent 9\quad  Appendix: Reciprocity for Kleinian groups
 
\section{Introduction}
 
Let $\cM_{g,n}$ be the moduli space of Riemann surfaces
of genus $g$ with $n$ punctures.

From a complex perspective, moduli space is hyperbolic.
For example, $\cM_{g,n}$ is abundantly populated by
immersed holomorphic disks of constant curvature $-1$ in the
Teichm\"uller (=Kobayashi) metric.

When $r=\dim_{\cx} \cM_{g,n}$ is greater than one,
however, $\cM_{g,n}$ carries no complete metric of bounded negative
curvature.
Instead, Dehn twists give chains of subgroups
$\zed^r \subset \pi_1(\cM_{g,n})$ reminiscent of
flats in symmetric spaces of rank $r>1$.

In this paper we introduce a new K\"ahler metric 
on moduli space that
exhibits its hyperbolic tendencies
in a form compatible with higher rank.

\demo{Definitions}
Let $(M,g)$ be a K\"ahler manifold.
An $n$-form $\alpha$ is  $d$(\/{\it bounded}\/)  if
$\alpha=d\beta$ for some bounded $(n-1)$-form $\beta$.
The space $(M,g)$ is {\it K{\rm \"{\it a}}hler hyperbolic} if:
\begin{itemize}
        \item[1.]
On the universal cover $\widetilde{M}$,
the K\"ahler form $\omega$ of the pulled-back metric 
$\widetilde{g}$ is $d$(bounded);
        \item[2.]
$(M,g)$ is complete and of finite volume;
        \item[3.]
The sectional curvature of $(M,g)$ is bounded above
        and below; and
        \item[4.]
The injectivity radius of $(\widetilde{M},\widetilde{g})$ is bounded below.
\end{itemize}
Note that (2--4) are automatic if $M$ is compact. 
\enddemo

The notion of a K\"ahler hyperbolic manifold was introduced 
by Gromov.  
Examples include compact K\"ahler manifolds of negative
curvature, products of such manifolds,
and finite volume quotients of Hermitian symmetric
spaces with no compact or Euclidean 
factors \cite{Gr}.

In this paper we show: 

\proclaimtitle{K\"ahler hyperbolic}
\proclaim{Theorem} 
\label{thm:kh}
The Teichm{\rm \"{\it u}}ller metric on moduli space is comparable to
a K{\rm \"{\it a}}hler metric $h$ such that
$(\cM_{g,n},h)$ is K{\rm \"{\it a}}hler hyperbolic{\rm .}
\endproclaim 

\demo{The bass note of Teichm{\rm \"{\it u}}ller space}
The universal cover of $\cM_{g,n}$ is the
{\it Teichm\"uller space} $\cT_{g,n}$.
Recall that the Teichm\"uller metric gives
{\it norms} $\|\cdot\|_T$ on the tangent and cotangent bundles to
$\cT_{g,n}$.
The analogue of the lowest eigenvalue of the Laplacian for such 
a metric is:
$$
        \lambda_0(\cT_{g,n}) \EQ
        \inf_{f \mem C^\infty_0(\cT_{g,n})}
        \left.
        \int \|df\|_T^2 \,dV \right/ \int |f|^2 \,dV,
$$
where $dV$ is the volume element of unit norm.
\enddemo

\proclaim{{C}orollary} We have
$\lambda_0(\cT_{g,n}) > 0$ in the Teichm{\rm \"{\it u}}ller metric{\rm .}
\endproclaim 

\demo{Proof}
The K\"ahler metric $h$ is comparable to the Teichm\"uller metric, 
so it suffices to bound $\lambda_0(\cT_{g,n},h)$.
Since the K\"ahler form $\omega$ for $h$ is
$d$(bounded), say $\omega=d\theta$, the volume form
$\omega^n=d\eta = d(\theta \wedge \omega^{n-1})$ is also
$d$(bounded).   Using the Cauchy-Schwarz inequality we 
then obtain
\begin{eqnarray*}
        \brackets{f,f} 
        & = & 
        \int f^2 \, \omega^n
        = \int f^2 \,d\eta
        = - \int 2f \, df \wedge \eta \\
        & \le &
        C \brackets{f,f}^{1/2} \brackets{df,df}^{1/2} .
\end{eqnarray*}
The lower bound
$\brackets{df,df}/\brackets{f,f} \GE 1/C^2 > 0$ follows,
yielding $\lambda_0>0$.
\enddemo 

\proclaimtitle{Complex isoperimetric inequality}
\proclaim{{C}orollary}
For any compact complex submanifold $N^{2k} \subset \cT_{g,n}${\rm ,} we have
$$
        \vol_{2k}(N) \LE C_{g,n} \cdot \vol_{2k-1}(\bdry N)
$$
in the Teichm{\rm \"{\it u}}ller metric.
\endproclaim 

\demo{Proof}
Passing to the equivalent K\"ahler hyperbolic metric $h$, 
Stokes' theorem yields:
$$
        \vol_{2k}(N) = \int_N \omega^k
                = \int_{\bdry N} \theta \wedge \omega^{k-1}
                = O(\vol_{2k-1}(\bdry N)),
$$
since $\theta \wedge \omega^{k-1}$ is a bounded $2k-1$ form.
\enddemo

\noindent (These two corollaries also hold in the
Weil-Petersson metric, since its K\"ahler form is
$d$(bounded) by Theorem \ref{thm:prim} below.)

\demo{The Euler characteristic}
Gromov shows the Laplacian on the universal
        cover $\tilde{M}$ of a K\"ahler hyperbolic
        manifold $M$ is positive on $p$-forms,
        so long as $p \neq n=\dim_{\cx} M$.
The $L^2$-cohomology of $\widetilde{M}$ is
therefore concentrated in the middle dimension $n$.
Atiyah's $L^2$-index formula for the Euler characteristic
(generalized to complete manifolds of finite volume
and bounded geometry by Cheeger and 
Gromov \cite{CG}) then yields
$$
        \sign \chi(M^{2n}) = (-1)^n   .
$$
In particular, Chern's conjecture on the sign of
$\chi(M)$ for closed negatively curved manifolds holds
in the K\"ahler setting.  
See  \cite[\S 2.5A]{Gr}.

For moduli space we obtain: \enddemo

\proclaim{{C}orollary}
The orbifold Euler characteristic of moduli space satisfies 
$\chi(\cM_{g,n}) > 0$ if $\dim_{\cx} \cM_{g,n}$ is even{\rm ,}
and $\chi(\cM_{g,n}) < 0$ if $\dim_{\cx} \cM_{g,n}$ is odd{\rm .}
\endproclaim 

This corollary was previously known by explicit computations.
For example the Harer-Zagier formula gives
$$
        \chi(\cM_{g,1}) = \zeta(1-2g)
$$
for $g>2$, and this formula
alternates sign as $g$ increases \cite{HZ}.
\bigbreak
\begin{center}
 \BoxedEPSF{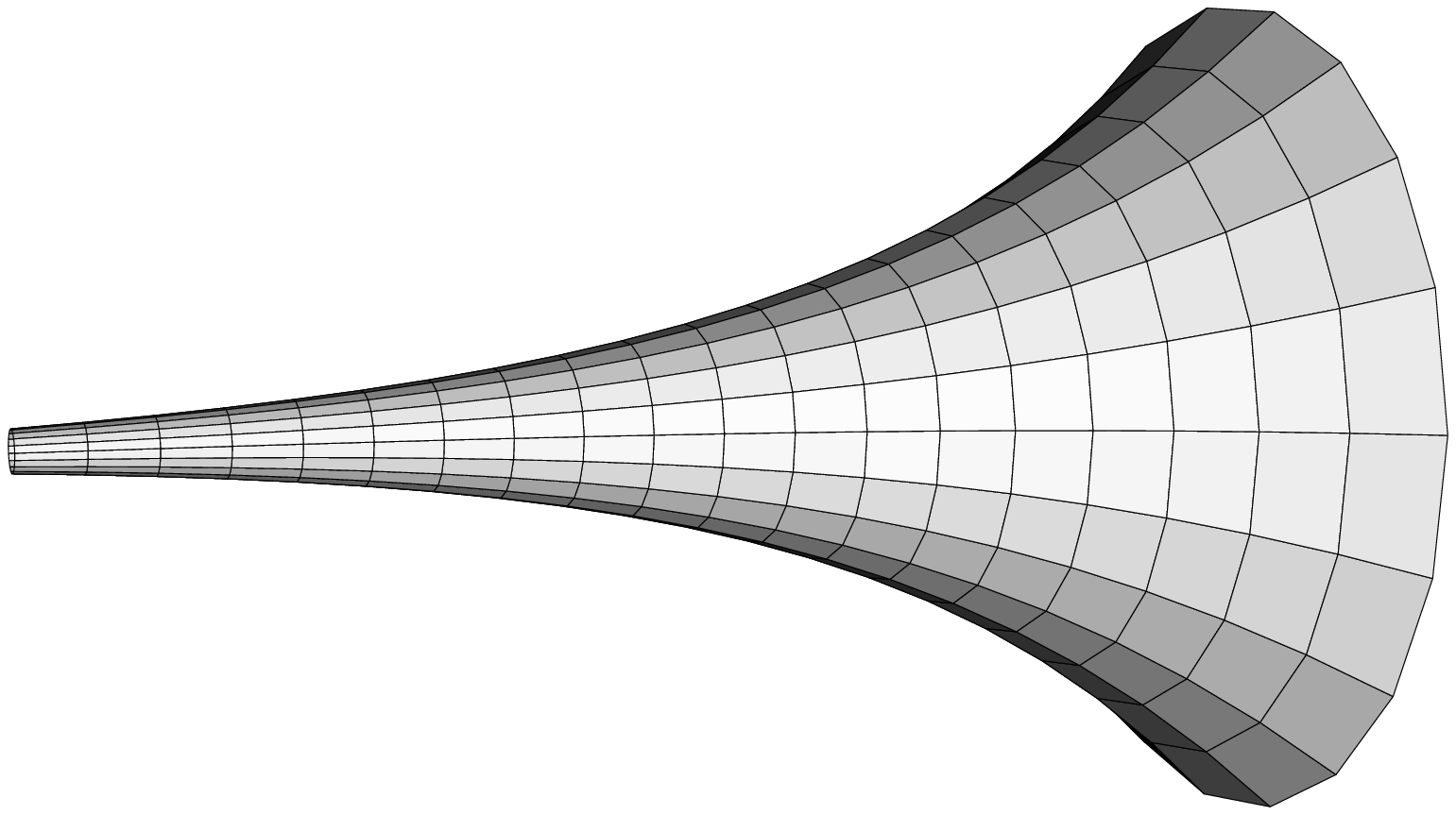 scaled 250}\hskip1in  \BoxedEPSF{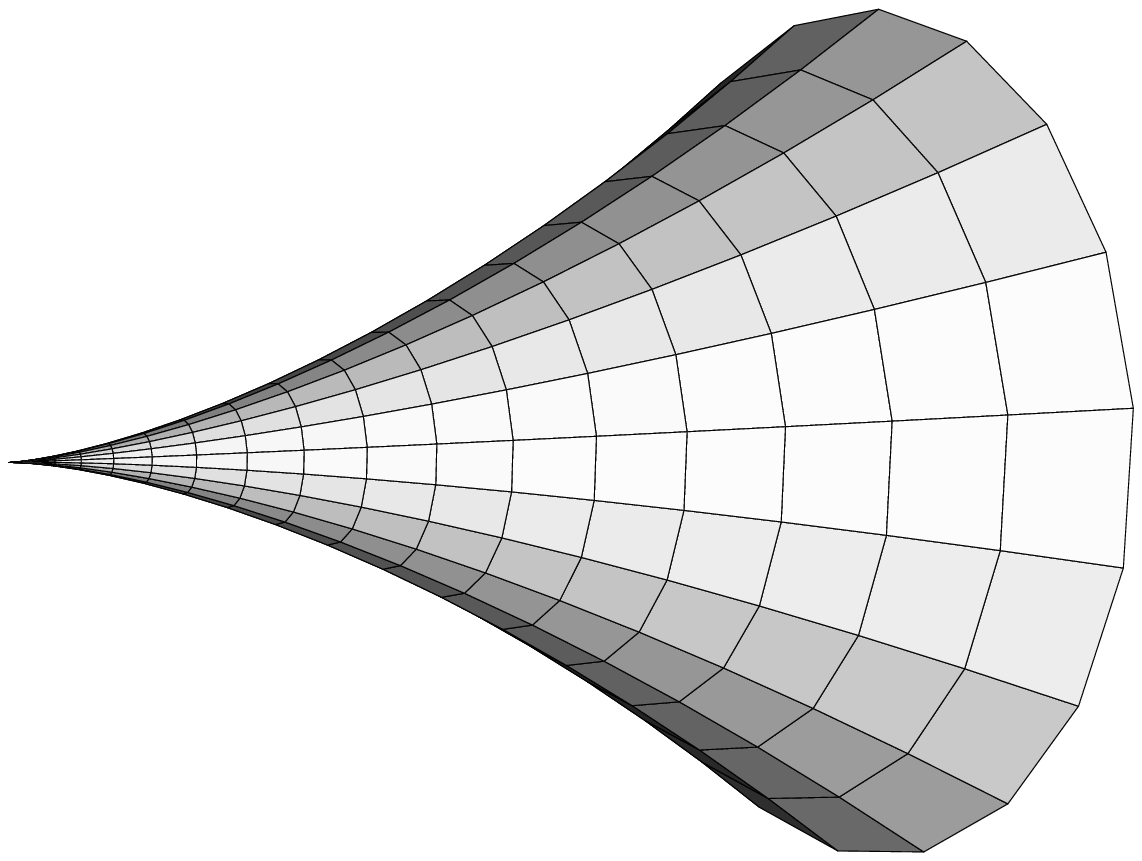 scaled 250}
\bigbreak
 {\ninepoint Figure 1. The cusp of moduli space in the Teichm\"uller and
Weil-Petersson metrics.}
\end{center}

\demo{Metrics on Teichm{\rm \"{\it u}}ller space}
To discuss the K\"ahler hyperbolic metric $h=\gol$
used to prove Theorem \ref{thm:kh},
we begin with the Weil-Petersson and Teichm\"uller metrics.  

Let $S$ be a hyperbolic Riemann surface   of genus $g$
with $n$ punctures,
and let $\Teich(S) \isom \cT_{g,n}$ be its
Teichm\"uller space.
The cotangent space $\T^{*}_X\Teich(S)$ is canonically
identified with the space $Q(X)$ of holomorphic quadratic
differentials $\phi(z) \,dz^2$ on $X \mem \Teich(S)$.
The Weil-Petersson and the Teichm\"uller metrics
correspond to the norms
\begin{eqnarray*}
        \|\phi\|_\WP^2 & = & \int_X \rho^{-2}(z) |\phi(z)|^2 \, |dz|^2
        \;\;\;\;{\rm  and}\\
        \|\phi\|_T & = & \int_X |\phi(z)| \, |dz|^2
\end{eqnarray*}
on $Q(X)$, where $\rho(z) |dz|$ is the hyperbolic metric on $X$.
The Weil-Petersson metric is K\"ahler, but the
Teichm\"uller metric is not even Riemannian when
$\dim_{\cx} \Teich(S) >1$.

To compare these metrics, consider the 
case of punctured tori with $\cT_{1,1} \isom \half \subset \cx$.
The Teichm\"uller metric on $\half$ is given by $|dz|/(2y)$,
while the Weil-Petersson metric is asymptotic to
$|dz|/y^{3/2}$ as $y \arrow \infty$.
Indeed, the Weil-Petersson 
symplectic form is given in
Fenchel-Nielsen length-twist coordinates by
$\omega_\WP = d\ell \wedge d\tau$,
and we have $\ell \asyto 1/y$ while $\tau \asyto x/y$.
Compare \cite{Mas}.

The cusp of the moduli space $\cM_{1,1} = \half/{\rm SL}_2(\zed)$ 
behaves like the surface of revolution
for $y=e^x$, $x<0$ in the Teichm\"uller metric; it is
complete and of constant negative curvature.
In Weil-Petersson geometry, on the other hand, 
the cusp behaves like the surface of revolution for $y=x^{3/2}$, $x>0$.
The Weil-Petersson metric on moduli space is
convex but incomplete, and its curvature
tends to $-\infty$ at the cusp.  See Figure~1.
\enddemo

\demo{A quasifuchsian primitive for the Weil-Petersson form}
Nevertheless the Weil-Petersson symplectic form 
$\omega_\WP$ is $d$(bounded),
and it serves as our point of departure for the
construction of a K\"ahler hyperbolic metric.
To describe a bounded primitive for $\omega_\WP$, recall
that the {\it Bers embedding} 
$$
        \beta_X : \Teich(\Sbar) \arrow Q(X)
        \isom \T^{*}_X \Teich(S) 
$$
sends Teichm\"uller space to a bounded domain in the
space of holomorphic quadratic differentials on $X$ (\S\ref{sec:basics}).
\enddemo

\proclaim{Theorem}
\label{thm:prim}
For any fixed $Y \mem \Teich(\Sbar)${\rm ,} the $1$\/{\rm -}\/form 
$$
        \theta_\WP(X) = -\beta_X(Y)
$$
is bounded in the Teichm{\rm \"{\it u}}ller and Weil\/{\rm -}\/Petersson metrics{\rm ,}
and satisfies $d(i\theta_\WP) = \omega_\WP${\rm .}
\endproclaim

The complex projective structures on $X$ are an affine
space modeled on $Q(X)$, and we can also write
$$
        \theta_\WP(X) = \sigf(X) - \sigqf(X,Y),
$$
where $\sigf(X)$ and $\sigqf(X,Y)$ are the Fuchsian and
quasifuchsian projective structures on $X$
(the latter coming from Bers' simultaneous uniformization
of $X$ and $Y$).
The 1-form $\theta_\WP$ is bounded by
Nehari's estimate for the Schwarzian derivative of
a univalent map (\S\ref{sec:dbounded}). 

Theorem \ref{thm:prim} is inspired by the formula
\begin{equation}
\label{eq:TZ}
        d(\sigma_F(X) - \sigma_S(X)) = -i\omega_\WP
\end{equation}
discovered by Takhtajan and Zograf, 
where the projective structure $\sigma_S(X)$
comes from a Schottky uniformization of $X$
\cite[Thm. 3]{Tak}, \cite{TZ};
see also \cite{Iv1}.
The proof of (\ref{eq:TZ}) by Takhtajan and Zograf
leads to remarkable results on the classical problem
of accessory parameters.
It is based on an
explicit K\"ahler potential for $\omega_\WP$ coming from
the {\it Liouville action} in string theory.
Unfortunately Schottky uniformization makes the 
1-form $\sigma_F(X) - \sigma_S(X)$ unbounded.

Our proof of Theorem \ref{thm:prim} is quite different and
invokes a new duality for Bers embeddings
which we call {\it quasifuchsian reciprocity} (\S\ref{sec:reciprocity}).    

\proclaim{Theorem}
Given $(X,Y) \mem \Teich(S) \times \Teich(\Sbar)${\rm ,}
the derivatives of the Bers embeddings
\begin{eqnarray*}
        D\beta_X  :  \T_Y \Teich(\Sbar) & \arrow & \T_X^*\Teich(S)
        \;\;\;\;\;{\rm and}\\
        D\beta_Y  :  \T_X \Teich(S) & \arrow & \T_Y^*\Teich(\Sbar)
\end{eqnarray*}
are adjoint linear operators\/{\rm ;} that is{\rm ,}
$D\beta_X^* = D\beta_Y$.
\endproclaim

Using this duality{\rm ,} we find that $d\theta_\WP(X)$ is
independent of the choice of~$Y${\rm .}
Theorem {\rm \ref{thm:prim}} then follows easily by setting $Y=\Xbar${\rm .}

In the Appendix we formulate a reciprocity law for
general Kleinian groups, and sketch a new proof of
the Takhtajan-Zograf formula (\ref{eq:TZ}).

\demo{The $1/\ell$ metric}
For any closed geodesic $\gamma$ on $S$, let
$\ell_\gamma(X)$ denote the length of the corresponding
hyperbolic geodesic on $X \mem \Teich(S)$.
A sequence $X_n \mem \cM(S)$ tends to infinity if and only
if $\inf_\gamma \ell_\gamma(X_n) \arrow 0$
\cite{Mum}.
This behavior motivates our use of the reciprocal length functions 
$1/\ell_\gamma$ to define a complete K\"ahler metric 
$\gol$ on moduli space.

To  begin the definition, let $\Log : \reals_+ \arrow [0,\infty)$ be a 
$C^\infty$ function such that
$$
        \Log(x) = 
                \left\{ \begin{array}{ll} \log(x) &\hbox{if $x \ge 2$,}\\
                              0       &\hbox{if $x \le 1$}.
                \end{array}\right.
$$
The {\it $1/\ell$ metric} $\gol$ is then defined,
for suitable small $\varepsilon$ and $\delta$, 
by its K\"ahler form
\begin{equation}
\label{eq:intro:wol}
                \wol \EQ
                w_\WP -i\delta
                \sumg
                \del \delbar \Log \frac{\varepsilon}{\ell_\gamma}
                \cdot
\end{equation}
The sum above is over primitive short geodesics $\gamma$
on $X$; at most $3|\chi(S)|/2$ terms occur in the sum.

Since $\gol$ is obtained by modifying the
Weil-Petersson metric, it is useful to have a
comparison between $\|v\|_T$ and $\|v\|_\WP$
based on short geodesics. \enddemo

\proclaim{Theorem}
For all $\varepsilon > 0$ sufficiently small{\rm ,} we have\/{\rm :}\/
\begin{equation}
\label{eq:intro:compare}
        \|v\|_T^2 \;\bab\;
        \: \|v\|_\WP^2 \:+\:
        \sumg |(\del \log \ell_\gamma)(v)|^2  .
\end{equation}
\endproclaim

This estimate {\rm (\S\ref{sec:comparison})}
is based on a thick\/{\rm -}\/thin decomposition for
quadratic differentials {\rm (\S\ref{sec:thick:thin}).}

\demo{Proof of Theorem {\rm \ref{thm:kh}}}
We can now outline the proof that $h=\gol$ is
K\"ahler hyperbolic and
comparable to the Teichm\"uller metric.

We begin by showing that any geodesic length function
is almost pluriharmonic (\S\ref{sec:ah}); more precisely,
$$
        \|\del \delbar (1/\ell_\gamma)\|_T = O(1).
$$
This means the term $\del \delbar \Log (\varepsilon/\ell_\gamma)$ 
in the definition (\ref{eq:intro:wol}) of $\wol$
can be replaced by
$(\del \Log \ell_\gamma) \wedge (\delbar \Log \ell_\gamma)$
with small error.
Using the relation between the Weil-Petersson and
Teichm\"uller metrics given by (\ref{eq:intro:compare}),
we then obtain
the comparability estimate $\gol(v,v) \bab \|v\|_T^2$.
This estimate implies moduli space is complete and of finite volume in
the metric $\gol$, because the same statements hold
for the Teichm\"uller metric.

To show $\wol$ is $d$(bounded), we note that
$d(i\tol) = \wol$ where
$$
        \tol \EQ
        \theta_\WP - \delta 
                \sumg
                \del \Log \frac{\varepsilon}{\ell_\gamma} \cdot
$$
The first term $\theta_\WP$ is bounded by Theorem \ref{thm:prim},
and the remaining terms are bounded by
basic estimates for the gradient of geodesic length.

Finally we observe that $\ell_\gamma$ and
$\theta_\WP$ can be
extended to holomorphic functions on
the complexification of $\Teich(S)$.  
Local uniform bounds on these holomorphic functions control all
their derivatives, and yield the desired bounds on the
curvature and injectivity radius of $\gol$
(\S\ref{sec:volume}).
\enddemo

\demo{The $1/d$ metric and domains in the plane}
To conclude we mention a parallel discussion of 
a K\"ahler metric $\god$ comparable to the
hyperbolic metric $g_H$ on a bounded domain
$\Omega \subset \cx$ with smooth boundary.

The (incomplete) 
Euclidean metric $g_E$ on $\Omega$ is defined by the K\"ahler form
$$
        \omega_E 
        \EQ \frac{i}{2} \, dz \wedge \dzbar .
$$
A well-known argument (based on the Koebe $1/4$-theorem)
gives for $v \mem \T_z\Omega$ the estimate
\begin{equation}
\label{eq:od}
        \|v\|_H^2 \;\bab\; \frac{\|v\|^2_E}{d(z,\bdry \Omega)^2},
\end{equation}
where $d(z,\bdry \Omega)$ is the Euclidean distance
to the boundary \cite{BP}.

Now consider the {\it $1/d$ metric} $\god$, defined for small
$\varepsilon$ and $\delta$ by the K\"ahler form
$$
        \wod(z) = \omega_E(z) \:+\:
        i\,\delta \, \del \delbar 
                \Log \frac{\varepsilon}{d(z,\bdry \Omega)} \cdot
$$
We claim that for suitable $\varepsilon$ and $\delta$,
the metric $\god$ is comparable to
the hyperbolic metric $g_H$.
\enddemo

\demo{Sketch of the proof}
Since $\bdry \Omega$ is smooth,
the function $d(z)=d(z,\bdry \Omega)$ is also smooth near the boundary
and satisfies 
$\|\del \delbar d\|_H = O(d^2)$.
Thus for $\varepsilon>0$ sufficiently small,
$\del \delbar \Log(\varepsilon/d)$ is dominated by the
gradient term $(\del d \wedge \delbar d)/d^2$.
Since $|(\del d)(v)|$ is comparable to the Euclidean
length $\|v\|_E$, by (\ref{eq:od}) we find
$g_H \bab \god$.
\hfill\qed\enddemo

Like the function $1/d(z,\bdry \Omega)$, the
reciprocal length functions
$1/\ell_\gamma(X)$ measure
the distance from $X$ to the boundary of moduli space,
rendering the metric\linebreak 
$\gol$ complete and comparable to
the Teichm\"uller (=Kobayashi) metric on~$\cM(S)$.

\demo{References}
The curvature and convexity of the 
Weil-Petersson metric and the behavior of
geodesic length-functions are discussed in \cite{Wol1}
and \cite{Wol2}.
For more on $\pi_1(\cM_{g,n})$, its subgroups
and parallels with lattices in Lie groups,
see \cite{Iv2}, \cite{Iv3}.
The hyperconvexity of Teichm\"uller space, which is
related to K\"ahler hyperbolicity, is
established by Krushkal in \cite{Kru}.
\enddemo
 
{\it Acknowledgements}.
I would like to thank Gromov for posing the question of the
K\"ahler hyperbolicity of moduli space,
and Takhtajan for explaining his work with Zograf several years ago.
Takhtajan also provided useful and insightful 
remarks when this paper was first circulated,
leading to the Appendix.  
\medbreak

{\it Notation}.
We use the standard notation $A=O(B)$ to mean
$A \le CB$, and $A \bab B$ to mean
$A/C < B < CA$, for some constant $C>0$.
Throughout the exposition, the constant $C$ is allowed to
depend on $S$ but it is otherwise universal.
In particular, all bounds will be uniform over
the entire Teichm\"uller space of $S$ unless
otherwise stated. 

\section{Teichm\"uller space}
\label{sec:basics}

This section reviews basic definitions and constructions
in Teichm\"uller theory;
for further background see
\cite{Gd}, 
\cite{IT},
\cite{Le}, and
\cite{Nag}.

\demo{The hyperbolic metric} 
A Riemann surface $X$ is {\it hyperbolic}  if it is covered
by the upper halfplane $\half$.
In this case the metric
$$
        \rho = \frac{|dz|}{\Im z}
$$
on $\half$ descends to the {\it hyperbolic metric} on $X$,
a complete metric of constant curvature $-1$. \enddemo

\demo{The Teichm{\rm \"{\it u}}ller metric}
Let $S$ be a hyperbolic Riemann surface.
A Riemann surface $X$ is {\it marked} by $S$ if it is
equipped with a quasiconformal
homeomorphism $f : S \arrow X$.
The {\it Teichm\"uller metric} on marked surfaces is defined~by
$$
        d((f : S \arrow X),(g : S \arrow Y)) = \frac{1}{2} \inf \log K(h),
$$
where $h : X \arrow Y$ ranges over all quasiconformal maps
isotopic to $g \compos f^{-1}$ {\rm rel} ideal boundary,
and $K(h) \ge 1$ is the dilatation of $h$.
Two marked surfaces are {\it equivalent} if their Teichm\"uller distance
is zero; then there is a conformal map
$h : X \arrow Y$ respecting the markings.  The metric
space of equivalence classes is the {\it Teichm\"uller space}
of $S$, denoted $\Teich(S)$.

Teichm\"uller space is naturally a complex manifold.
To describe its tangent and cotangent spaces, let
$Q(X)$ denote the Banach space
of holomorphic quadratic differentials $\phi = \phi(z) \,dz^2$
on $X$ for which the $L^1$-norm
$$
        \|\phi\|_T \EQ \int_X | \phi |
$$
is finite; and  let $M(X)$ be the space of 
$L^\infty$ measurable Beltrami differentials $\mu(z) \,\dzbar/dz$
on $X$.  
There is a natural pairing between $Q(X)$ and $M(X)$ 
given~by
$$
        \brackets{\phi,\mu} = \int_X \phi(z) \mu(z) \, dz \,\dzbar.
$$
A vector $v \mem \T_X \Teich(S)$ is represented by a
Beltrami differential $\mu \mem M(X)$, and its
{\it Teichm\"uller norm} is given by 
$$
        \|\mu\|_T = \sup \{\Re \brackets{\phi,\mu} \st \|\phi\|_T = 1\} .
$$
We have the isomorphism:
$$
        \T_X\Teich(S) \isom Q(X)^* \isom M(X)/Q(X)^\perp ,
$$
and $\|\mu\|_T$ gives the infinitesimal form of the Teichm\"uller metric.
\enddemo

\demo{Projective structures}
A complex {\it projective structure} on $X$ is a
subatlas of charts whose transition functions
are M\"obius transformations.
The space of projective surfaces marked by $S$
is naturally a complex manifold $\Proj(S) \arrow \Teich(S)$
fibering over Teichm\"uller space.
The Fuchsian uniformization, $X=\half/\Gamma(X)$, determines
a canonical section
$$
        \sigf : \Teich(S) \arrow \Proj(S).
$$
This section is real analytic but {\it not} holomorphic.

Let $P(X)$ be the Banach space of holomorphic quadratic 
differentials on $X$ with finite $L^\infty$-norm
$$
        \|\phi\|_\infty = \sup_X \rho^{-2}(z) |\phi(z)| .
$$
The fiber $\Proj_X(S)$ of $\Proj(S)$ over $X \mem \Teich(S)$ is an
affine space modeled on $P(X)$.
That is, given $X_0 \mem \Proj_X(S)$ and $\phi \mem P(X)$,
there is a unique $X_1 \mem \Proj_X(S)$ and a 
conformal map $f : X_0 \arrow X_1$ respecting markings, such that
$Sf=\phi$.  Here $Sf$ is the {\it Schwarzian derivative}
$$
        Sf(z) = 
        \left(\frac{f''(z)}{f'(z)} \right)'
        - \frac{1}{2} \left(\frac{f''(z)}{f'(z)} \right)^2 \, dz^2 .
$$
Writing $X_1 = X_0 + \phi$, we have $\Proj_X(S) = \sigf(X) + P(X)$.
\enddemo

\demo{Nehari{\rm '}s bound}
A {\it univalent function} is an injective, holomorphic map 
$f : \half \arrow \chat$.
The bounds of the next result \cite[\S 5.4]{Gd}
play a key role in proving
universal bounds on  the geometry of $\Teich(S)$. \enddemo

\proclaimtitle{Nehari}
\proclaim{Theorem} 
\label{thm:Nehari}
Let $Sf$ be the Schwarzian derivative of a holomorphic map
$f : \half \arrow \chat$.  Then we have the implications\/{\rm :}
$$
        \|Sf\|_\infty < 1/2 \;\implies\;
        \hbox{ {\rm (}$f$ is univalent\/{\rm )}}
        \;\implies\;
        \|Sf\|_\infty < 3/2.
$$
\endproclaim

\demo{Quasifuchsian groups}
The space $QF(S)$ of marked quasifuchsian groups provides
a complexification of $\Teich(S)$ that plays a crucial role
in the sequel.

Let $\chat = \half \cup \lowerhalf \cup  \Rinf$
denote the partition of the Riemann sphere into the
upper and lower halfplanes and the circle
$\Rinf=\reals \cup \{\infty\}$ .
Let $S = \half/\Gamma(S)$ be a presentation of $S$ as the
quotient $\half$ by the action of a Fuchsian group 
$\Gamma(S) \subset {\rm PSL}_2(\reals)$.

Let $\Sbar = \lowerhalf/\Gamma$ denote the complex conjugate of $S$.
Any Riemann surface $X \mem \Teich(S)$ also has
a complex conjugate $\Xbar \mem \Teich(\Sbar)$,
admitting an anticonformal map $\Xbar \arrow X$
compatible with marking.

The {\it quasifuchsian space} of $S$ is defined by
$$
        QF(S) = \Teich(S) \times \Teich(\Sbar) .
$$
The map $X \mapsto (X,\Xbar)$ sends Teichm\"uller space
to the totally real {\it Fuchsian subspace}
$F(S) \subset QF(S)$, and thus $QF(S)$ is a complexification
of $\Teich(S)$.

The space $QF(S)$ parametrizes
marked quasifuchsian groups equivalent to $\Gamma(S)$, as follows.
Given  
$$
        (f : S \arrow X, g : \Sbar \arrow Y) \mem QF(S),
$$
we can pull back the complex structure from $X \cup Y$ to
$\half \cup \lowerhalf$, solve the Beltrami equation, and
obtain a quasiconformal map $\phi : \chat \arrow \chat$
such that:
\begin{itemize}
        \item
        $\phi$ transports the action of $\Gamma(S)$ 
        to the action of a Kleinian group 
        $\Gamma(X,Y) \subset {\rm PSL}_2(\cx)$; 
        \item
        $\phi$ maps $(\half \cup\lowerhalf,\Rinf)$ to
        $(\Omega(X,Y),\Lambda(X,Y))$, 
        where $\Lambda(X,Y)$ is a quasicircle; and
        \item
        there is an isomorphism
        $\Omega(X,Y)/\Gamma(X,Y) \isom X \cup Y$ such that
        $$
                \phi : (\half \cup \lowerhalf) \arrow \Omega(X,Y)
        $$
        is a lift of $(f \cup g) : (S \cup \Sbar) \arrow (X \cup Y)$.
\end{itemize}
Then $\Gamma(X,Y)$ is a {\it quasifuchsian group} equipped with a
conjugacy $\phi$ to $\Gamma(S)$.  Here $(X,Y)$ determines $\Gamma(X,Y)$ 
up to conjugacy in ${\rm PSL}_2(\cx)$, and $\phi$ up to
isotopy {\rm rel} $(\Rinf,\Lambda(X,Y)).$\footnote{
When $S$ has finite area, the limit set of $\Gamma(X,Y)$
coincides with $\Lambda(X,Y)$; in general it may be smaller.}

There is a natural holomorphic map
$$
        \sigma  :  \Teich(S) \times \Teich(\Sbar)
                \arrow \Proj(S) \times \Proj(\Sbar), 
$$
which records the projective structures on $X$ and $Y$
inherited from $\Omega(X,Y) \subset \chat$.
We denote the two coordinates of this map by
$$
        \sigma(X,Y) = (\sigqf(X,Y),\sigqfbar(X,Y)) .
$$

The {\it Bers embedding} $\beta_Y : \Teich(S) \arrow P(Y)$
is given by 
$$
        \beta_Y(X) \EQ \sigqfbar(X,Y)\,-\,\sigf(Y) .
$$
Writing $Y=\half/\Gamma(Y)$, we have
$\beta_Y(X) = Sf$, where $f : \half \arrow \Omega(X,Y)$ is
a Riemann mapping conjugating $\Gamma(Y)$ to $\Gamma(X,Y)$.
Amplifying Theorem \ref{thm:Nehari} we have: \enddemo

\proclaim{Theorem}
\label{thm:Bers}
The Bers embedding maps Teichm{\rm \"{\it u}}ller space to a bounded
domain in $P(Y)${\rm ,} with
$$
        B(0,1/2) \subset \beta_Y(\Teich(S)) \subset B(0,3/2),
$$
where $B(0,r)$ is the norm ball of radius $r$ in $P(Y)${\rm .}
The Teichm{\rm \"{\it u}}ller metric agrees with the Kobayashi metric
on the image of $\beta_Y${\rm .}
\endproclaim

See \cite[\S 5.4, \S 7.5]{Gd}.
(This reference has different constants, because there the hyperbolic 
metric $\rho$ is normalized to have curvature $-4$ instead of $-1$.)

\demo{Real and complex length}
Given a hyperbolic geodesic $\gamma$ on $S$, let
$\ell_\gamma(X)$ denote the hyperbolic length of the 
corresponding geodesic on $X \mem \Teich(S)$.
For $(X,Y) \mem QF(S)$, we can normalize coordinates on
$\chat$ so that the element
$g \mem \Gamma(X,Y)$ corresponding to $\gamma$
is given by $g(z)=\lambda z$, $|\lambda|>1$,
and so that $1$ and $\lambda$ belong to $\Lambda(X,Y)$.
By analytically continuing the logarithm from $1$ to
$\lambda$ along
$\Lambda(X,Y)$, starting with $\log(1)=0$, we obtain
the {\it complex length}
$$
        \cL_\gamma(X,Y) = \log \lambda = L+i\theta .
$$
In the hyperbolic 3-manifold $\half^3/\Gamma(X,Y)$,
$\gamma$ corresponds to a closed geodesic of length $L$ and
torsion $\theta$.

The group $\Gamma(X,Y)$ varies holomorphically
as a function of $(X,Y) \mem QF(S)$, so we have: \enddemo

\proclaim{Proposition}
The complex length $\cL_\gamma : QF(S) \arrow \cx$ is holomorphic{\rm ,}
and satisfies $\ell_\gamma(X)=\cL_\gamma(X,\Xbar)$.
\endproclaim 

\demo{The Weil\/{\rm -}\/Petersson metric}
Now suppose $S$ has finite hyperbolic area.
The {\it Weil-Petersson metric} is defined on the cotangent space
$Q(X) \isom T^*_X \Teich(S)$ by the $L^2$-norm
$$
        \|\phi\|_\WP^2 =
                \int_X \rho^{-2}(z) \, |\phi|^2 \, |dz|^2 .
$$
By duality we obtain a Riemannian metric
$g_\WP$ on the tangent space to $\Teich(S)$, and in fact 
$g_\WP$ is a K\"ahler metric. \enddemo

\proclaim{Proposition}
\label{prop:WP:T}
For any tangent vector $v$ to $\Teich(S)$ we have
$$
        \|v\|_\WP \le |2 \pi \chi(S)|^{1/2} \cdot \|v\|_T .
$$
\endproclaim 

\demo{Proof}
By Cauchy-Schwarz, if $\phi \mem Q(X)$ represents a cotangent
vector then we have
$$
        \|\phi\|_T =
        \int_X \frac{|\phi|}{\rho^2} \rho^2
        \le 
        \left( 
        \int_X 1 \cdot {\rho^2}
        \right)^{1/2}
        \left( 
        \int_X \frac{|\phi|^2}{\rho^4} \rho^2
        \right)^{1/2}
        = 
        |2\pi \chi(S)|^{1/2} \cdot \|\phi\|_\WP,
$$
where Gauss-Bonnet determines the hyperbolic area of $S$.
By duality the reverse inequality holds on the tangent space.
\enddemo

\section{$1/\ell$ is almost pluriharmonic}
\label{sec:ah}

In this section we begin a more detailed study of geodesic
length functions and prove a universal bound on
$\del \delbar (1/\ell_\gamma)$.

The Teichm\"uller metric $\|v\|_T$ on tangent vectors 
determines a norm $\|\theta\|_T$ for $n$-forms on $\Teich(S)$ by
$$
        \|\theta\|_T = \sup \{ |\theta(v_1,\ldots,v_n)| \st \|v_i\|_T = 1\},
$$
where the sup is over all $X \mem \Teich(S)$ and all 
$n$-tuples $(v_i)$ of unit tangent vectors
at $X$.

\proclaimtitle{Almost pluriharmonic}
\proclaim{Theorem}
\label{thm:ah}
Let $\ell_\gamma : \Teich(S) \arrow \reals_+$
be the length function of a closed geodesic on $S${\rm .}  Then
$$
        \|\del \delbar (1/\ell_\gamma)\|_T = O(1).
$$
The bound is independent of $\gamma$ and $S${\rm .}
\endproclaim

We begin by discussing the case where $S$ is an 
annulus and $\gamma$ is its core geodesic.
To simplify notation, set $\ell=\ell_\gamma$ and
$\cL=\cL_\gamma$.
Each annulus $X \mem \Teich(S)$ can be presented as a quotient:
$$
        X \EQ \half/\brackets{z \mapsto e^{\ell(X)} z}.
$$
The metric $|dz|/|z|$ makes $X$ into a right 
cylinder of area $A=\pi \ell$ and circumference $C=\ell(X)$;
the {\it modulus} of $X$ is the ratio
$$
        \mod(X) = \frac{A}{C^2} = \frac{\pi}{\ell(X)} \cdot
$$

Given a pair of Riemann surfaces 
$(X,Y) \mem \Teich(S) \times \Teich(\Sbar)$ we can glue $X$ to $Y$
along their ideal boundaries (which are canonically identified
using the markings by $S$) to obtain a complex torus
$$
        T(X,Y) \EQ
        X \cup (\bdry X = \bdry Y) \cup Y  
        \isom
        \cx^*/\brackets{e^{\cL(X,Y)}} ,
$$
where $\cL(X,Y)$ is the complex length introduced in
Section~\ref{sec:basics}. 
This torus is simply the quotient Riemann surface for the
Kleinian group 
$$
        \Gamma(X,Y) \isom \brackets{z \mapsto e^{\cL(X,Y)}z}.
$$

The metric $|dz|/|z|$ makes $T(X,Y)$ into
a flat torus with area $A=2\pi \Re \cL$ in which
$\bdry X$ is represented by a geodesic loop
of length $C=|\cL|$.  
We define the {\it modulus} of the torus by
$$
        \mod(T(X,Y)) = \frac{A}{C^2}
                = \Re \frac{2\pi}{\cL(X,Y)}  \cdot
$$
Note that $T(X,\Xbar)$ is obtained by doubling the
annulus $X$, and $\mod(T(X,\Xbar)) = 2\mod(X)$.

\proclaim{Lemma}  
\label{lem:mod}
If the Teichm{\rm \"{\it u}}ller distance from $X$ to $\Ybar$
is bounded by $1${\rm ,} then
$$
        \mod(T(X,Y)) = \mod(X) + \mod(Y) + O(1) .
$$
\endproclaim 

\demo{Proof}
Since $d_T(X,\Ybar) \le 1$, there is a $K$-quasiconformal
map from $T(X,\Xbar)$ to $T(X,Y)$ with $K=O(1)$.
The annuli $X, Y \subset T(X,Y)$ are thus
separated by a pair of $K$-quasicircles.
A quasicircle has bounded turning
\cite[\S 8.7]{LV}, with a bound controlled by $K$,
so we can find a pair of geodesic cylinders (with respect
to the flat metric on $T(X,Y)$) such that
$\bdry X = \bdry Y \subset A \cup B$ and $\mod(A)=\mod(B)=O(1)$; see Figure~2 .
(The cylinders $A$ and $B$ will be embedded if $\mod(X)$ and $\mod(Y)$
are large; otherwise they may be just immersed.)

\figin{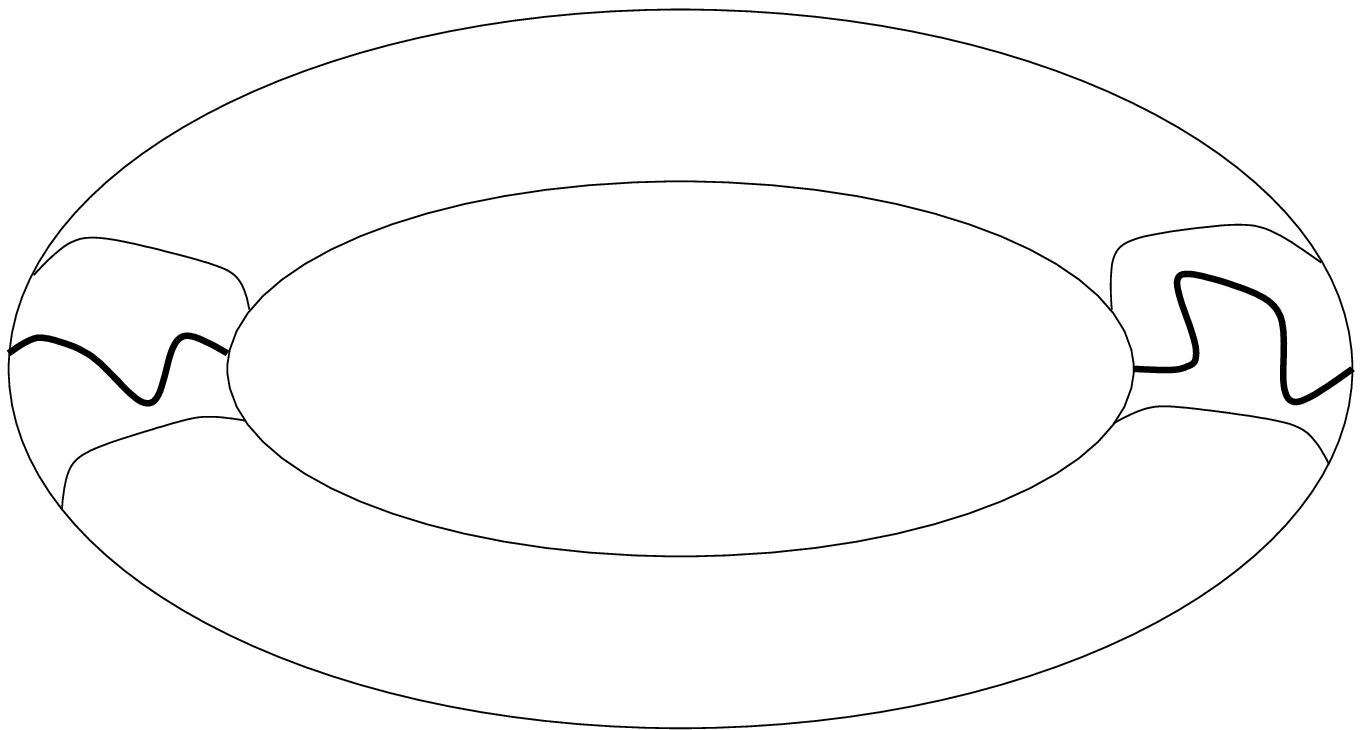}{500}
\vglue-2in 
\vglue.5pc
\phantom{blah}\hglue3.5in $T(X,Y)$
\vglue12pt
\vglue-2pc \phantom{more}\hskip2.25in $Y$
\vglue-1pc
\vglue.25in \phantom{give me room} \hglue.75in $A$\hglue3in $B$
\vglue-.25pc
\vglue.25in \phantom{blah} \hglue2.25in $X$
\vglue.25in
\centerline{\eightpoint Figure 2. Two annuli joined to form the torus $T(X,Y)$.}
\medbreak

The geodesic cylinders $X \cup A \cup B$ and $Y \cup A \cup B$
cover $T(X,Y)$ with bounded overlap, so their moduli sum
to $\mod(T(X,Y))+O(1)$.  Combining this fact with
monotonicity of the modulus \cite[\S 4.6]{LV},
we have
\begin{eqnarray*}
\mod(X) + \mod(Y) & \le & \mod(X \cup A \cup B) + \mod(Y \cup A \cup B) \\
        & =  & \mod(T) + O(1) .
\end{eqnarray*}
Similarly, we have
\begin{eqnarray*}
        \mod(T(X,Y)) & = &
        \mod(X-A-B) + \mod(Y-A-B) + O(1) \\
        & \le & \mod(X)+\mod(Y) + O(1),
\end{eqnarray*}
establishing the theorem.
\enddemo

{\it Proof of Theorem} \ref{thm:ah} (Almost pluriharmonic). 
We continue with the case of an annulus and its core geodesic
as above.  
Consider $X_0 \mem \Teich(S)$ and $v \mem \T_{X_0} \Teich(S)$ with
$\|v\|_T=1$.  Let $\Delta$ be the unit disk in $\cx$.
Using the Bers embedding of $\Teich(S)$ into $P(\Xbar_0)$
and Theorem \ref{thm:Bers},
we can find a holomorphic disk
$$
        \iota : (\Delta,0) \arrow (\Teich(S),X_0),
$$
tangent to $v$ at the origin, such that the Teichm\"uller
and Euclidean metrics are comparable on $\Delta$,
and $\diam_T(\iota(\Delta)) \le 1$.
(For example, we can take $\iota(s)=sv/10$ using the
linear structure on $P(\Xbar_0)$.)

Let $X_s = \iota(s)$ and
$Y_t = \Xbar_\tbar \mem \Teich(\Sbar)$;
then $(X_s,Y_t) \mem QF(S)$ is a holomorphic function of 
$(s,t) \mem \Delta^2$.
Set 
$$
        M(X,Y)=\mod(T(X,Y)) = \Re \frac{2\pi}{\cL(X,Y)},
$$
and define
$f : \Delta^2 \arrow \reals$ by
$$
        f(s,t) = M(X_s,Y_t) - M(X_s,Y_0) - 
                M(X_0,Y_t) + M(X_0,Y_0) .
$$
By Lemma \ref{lem:mod} above,   $f(s,t)=O(1)$.
On the other hand, $\cL(X,Y)$ is
holomorphic, so $f(s,t)$ is pluriharmonic.
Thus the bound $f(s,t)=O(1)$ controls the full 2-jet
of $f(s,t)$ at $(0,0)$; in particular,
$$
        \left.
        \frac{\del^2 f(s,t)}{\del s \,\del t} \right|_{0,0}
        = O(1) .
$$
By letting $g(s)=f(s,\sbar)$, it follows that
$(\del \delbar g)(0) = O(1)$
in the Euclidean metric on $\Delta$.
On the other hand,
$$
        \del \delbar g(s) = \del \delbar M(X_s,\Xbar_s) =
        \del \delbar (\pi/\ell(X_s)),
$$
since the remaining terms in the expression for $f(s,\sbar)$ are
pluriharmonic in $s$.  Thus
$\|\del \delbar (1/\ell)\|_T = O(1)$, and the proof is
complete for annuli.

To treat the case of general $(S,\gamma)$,
let $\widetilde{S} \arrow S$ be the annular covering space determined by
$\brackets{\gamma} \subset \pi_1(S)$, and let
$\pi : \Teich(S) \arrow \Teich(\widetilde{S})$ be the
holomorphic map obtained by lifting complex structures.
Then we have:
$$
        \|\del \delbar (1/\ell_\gamma)\|_T
        =
        \|\pi^*(\del \delbar (1/\ell))\|_T
        \le
        \|\del \delbar (1/\ell)\|_T
        = O(1),
$$
since holomorphic maps do not expand the 
Teichm\"uller ($=$Kobayashi) metric.
\phantom{timeforagame}\hfill\qed

\demo{{R}emark}
It is known that on finite-dimensional Teichm\"uller
spaces, $\ell_\gamma$ is
strictly plurisubharmonic \cite{Wol2}. \enddemo

\section{Thick-thin decomposition of quadratic differentials}
\label{sec:thick:thin}  

Let $S$ be a hyperbolic surface of finite area,
and let $\phi \mem Q(X)$ be a quadratic differential on
$X \mem \Teich(S)$.
In this section we will present a canonical decomposition
of $\phi$ adapted to the short geodesics $\gamma$ on $X$.

To each $\gamma$ we will associate a 
residue $\Res_\gamma : Q(X) \arrow \cx$
and a differential
$\phi_\gamma \mem Q(X)$ proportional to
$\del \log \ell_\gamma$ with
$\Res_\gamma(\phi_\gamma) \approx 1$.  We will then show:

\proclaimtitle{Thick-thin}
\proclaim{Theorem} 
\label{thm:thick:thin}
For $\varepsilon>0$ sufficiently small{\rm ,}
any $\phi \mem Q(X)$ can be uniquely expressed in the form
\begin{equation}
\label{eq:res}
        \phi = \phi_0  + \sumg
                a_\gamma \phi_\gamma 
\end{equation}
with $\Res_\gamma(\phi_0)=0$ for all $\gamma$ in the sum
above{\rm .}
Each term $\phi_0$ and $a_\gamma \phi_\gamma$ 
has Teichm{\rm \"{\it u}}ller norm $O(\|\phi\|_T)${\rm .}
\endproclaim

We will also show that 
$\|\phi_0\|_\WP \bab \|\phi_0\|_T$ (Theorem~\ref{thm:Ceps}).
Thus the thick-thin decomposition accounts for
the discrepancy between the Teichm\"uller and
Weil-Petersson norms on $Q(X)$ in terms of short geodesics on $X$.

\demo{The quadratic differential $\del \log \ell_\gamma$}
Let $\gamma$ be a closed hyperbolic geodesic on $S$.
Given $X \mem \Teich(S)$, let $\pi : X_\gamma \arrow X$
be the covering space
corresponding to $\brackets{\gamma} \subset \pi_1(S)$.
We may identify $X_\gamma$ with a round annulus
$$
        X_\gamma \isom A(R) =
                \{z \st R^{-1} < |z| < R\}.
$$
By requiring that $\widetilde{\gamma} \subset X_\gamma$ and $S^1$ 
agree as oriented loops,
we can make this identification unique up to rotations.

Consider the natural 1-form $\theta_\gamma = dz/z$  on
$X_\gamma$.
In the $|\theta|$-metric, $X_\gamma$ is a right cylinder of
circumference $C=2\pi$ and area $A=4\pi \log R$.
Thus we have
\begin{eqnarray*}
        \mod(X_\gamma)  & = &
        \frac{A}{C^2} 
        \EQ
        \frac{\log R}{\pi} 
        \EQ 
        \frac{\pi}{\ell_\gamma(X)} , \;\;\;{\rm and} \\
        \|\theta^2\|_T & = & A 
        \EQ 
        \frac{4\pi^3}{\ell_\gamma(X)} \cdot
\end{eqnarray*}

Define $\phi_\gamma \mem Q(X)$ by
$$
        \phi_\gamma = \pi_*(\theta_\gamma^2) =
        \pi_* \left( \frac{dz^2}{z^2} \right) \cdot
$$
The importance of $\phi_\gamma$ comes from its well-known
connection to geodesic length:
\begin{equation}
\label{eq:dlogL}
        (\del \log \ell_\gamma)(X)
        \EQ
        -  \frac{\ell_\gamma(X)}{2\pi^3} 
        \, \phi_\gamma
\end{equation}
in $T^*_X \Teich(S) \isom Q(X)$ (cf. \cite[Thm.~3.1]{Wol2}).
\enddemo

\proclaim{Theorem}
\label{thm:ell}
The differential $(\del \log \ell_\gamma)(X)$ is proportional to
$\phi_\gamma${\rm .}  We have $\|\del \log \ell_\gamma\|_T \le 2${\rm ,} and
$\|\del \log \ell_\gamma\|_T \arrow 2$ as $\ell_\gamma \arrow 0${\rm .}
\endproclaim

\demo{Proof}
Equation (\ref{eq:dlogL}) gives the proportionality and
implies the bound
$$
        \|\del \log \ell_\gamma\|_T 
        =
        \frac{\ell_\gamma(X)}{2\pi^3} 
        \|\phi_\gamma\|_T \le
        \frac{\ell_\gamma(X)}{2\pi^3} \|\theta^2\|_T = 2.
$$
To analyze the behavior of  $\del \log \ell_\gamma$
when $\ell_\gamma(X)$ is small,
note that the collar lemma \cite{Bus}
provides a universal $\varepsilon_0>0$ such that for
$$
        T = \varepsilon_0 R,
$$
the map $\pi$ sends $A(T) \subset A(R)$ injectively into
a collar neighborhood of $\gamma$ on $X$.
Since $\int_{A(R)-A(T)} |\theta_\gamma|^2 = O(1)$,
we obtain
$$
        \|\phi_\gamma\|_T
        = \int_{\pi(A(T))} |\phi_\gamma| + O(1)
        = \|\theta^2\|_T + O(1) ,
$$
which implies $\|\del \log \ell_\gamma\|_T = 2 + O(\ell_\gamma)$.
\enddemo

\demo{The residue of a quadratic differential}
Let us define the {\it residue} of 
$\phi \mem Q(X)$ around $\gamma$ by
$$
        \Res_\gamma(\phi) = \frac{1}{2\pi i} 
                \int_{S^1} \frac{\pi^*(\phi)}{\theta_\gamma} \cdot
$$
In terms of the Laurent expansion
$$
        \pi^*(\phi) = \left( \sum_{-\infty}^\infty
                a_n z^n \right) \frac{dz^2}{z^2}
$$
 on $A(R)$, we have $\Res_\gamma(\phi)=a_0.$
\enddemo

\demo{Proof of Theorem {\rm \ref{thm:thick:thin} (Thick-thin)}}
To begin we will show that for
any $\gamma$ with $\ell_\gamma(X) < \epsilon$, we have
\begin{equation}
\label{eq:Res}
        \Res_\gamma(\phi) = 
        O \left( \frac{\|\phi\|_T}{\|\phi_\gamma\|_T} \right) .
\end{equation}
To see this, identify $X_\gamma$ with $A(R)$, set
$T=\epsilon_0 R$ as in the proof of Theorem \ref{thm:ell},
and consider the Beltrami coefficient on $A(T)$ given by 
$$
        \mu = 
        \frac{\overline{\theta_\gamma}^2}{|\theta_\gamma|}
        = \frac{z}{\zbar}\, \frac{d\zbar}{dz} \cdot
$$
Then we have:
\begin{equation}
\label{eq:AT}
        \Res_\gamma(\phi) = 
        \frac{1}{4\pi \log T}
        \int_{A(T)} \pi^*(\phi) \,\mu\,
        = 
        O\left(
        (\log T)^{-1} \int_{A(T)} |\pi^*\phi|
        \right). \hskip.4in
\end{equation}
Since $\pi|A(T)$ is injective, we have
$\int_{A(T)} |\pi^*\phi| = O(\|\phi\|_T)$
and $\log T \bab \|\phi_\gamma\|_T$, yielding (\ref{eq:Res}).

By similar reasoning, all
$\gamma$ and $\delta$ shorter than $\epsilon$ satisfy:
\begin{equation}
\label{eq:delta}
        \Res_{\gamma}(\phi_\delta) = 
        \left\{ \begin{array}{ll}
        1 & \hbox{if $\gamma=\delta$,}\\
        0 & \hbox{otherwise}
        \end{array}
        \right\} \;+\; 
        O \left( \frac{1}{\|\phi_\gamma\|} \right)
        \cdot
\end{equation}
Indeed, if $\delta \neq \gamma$ then
most of the mass of $|\phi_\delta|$ resides in the thin
part associated to $\delta$, which is disjoint from 
$\pi(A(T))$.  More precisely, we have
$\int_{A(T)} |\pi^*\phi_\delta| = O(1)$, and the desired
bound on $\Res_\gamma(\phi_\delta)$ follows from
(\ref{eq:AT}).  The estimate when $\delta=\gamma$
is similar, using the fact that
$\pi^*\phi_\gamma = \pi^*\pi_*(\theta_\gamma^2) \approx
\theta_\gamma^2$ on $A(T)$.

By (\ref{eq:delta}),
the matrix $\Res_\gamma(\phi_\delta)$ is close to the 
identity when $\epsilon$ is small, since
$\|\phi_\gamma\|_T^{-1} = O(\epsilon)$.
Therefore we have unique coefficients $a_\gamma$ satisfying
equation (\ref{eq:res}) in the statement of the theorem.

To estimate $|a_\gamma|$, we first use
the matrix equation
\begin{displaymath}
        [a_\gamma] = [\Res_\gamma \phi_\delta]^{-1}
        [\Res_\delta(\phi)]
\end{displaymath}
to obtain the bound
\begin{equation}
\label{eq:at}
        |a_\gamma| = O(\|\phi\|_T)
\end{equation}
from (\ref{eq:Res}).  
(Note that size of the matrix
$\Res_\gamma\phi_\delta$ is
controlled by the genus of $X$.)
Then we make the more precise estimate
\bigbreak
\centerline{${\displaystyle
        |a_\gamma|
        \bab 
        |\Res_\gamma(a_\gamma \phi_\gamma)|
        = 
        \left|
        \Res_\gamma(\phi) - \sum_{\delta \neq \gamma} 
                a_\delta \Res_\gamma(\phi_\delta)
        \right| 
        = O(\|\phi\|_T/\|\phi_\gamma\|_T)
}$}
\medbreak\noindent 
by (\ref{eq:Res}), (\ref{eq:delta}) and (\ref{eq:at}).
The bound $\|a_\gamma \phi_\gamma\|_T = O(\|\phi\|_T)$
follows.
  \enddemo

The bound on the terms in (\ref{eq:res}) above can be
improved when $\phi$ is also associated to a short geodesic.

\proclaim{Theorem}
\label{thm:gamma:medium}
If $\phi=\phi_\delta$ with
$2\varepsilon > \ell_\delta(X) > \varepsilon${\rm ,} then we have
$$
        \|a_\gamma \phi_\gamma\|_T
        \EQ O(\ell_\delta(X) \: \|\phi\|_T ) 
$$
in equation {\rm (\ref{eq:res}).}
\endproclaim

\demo{Proof}
For any $\gamma$ with $\ell_\gamma(X) < \varepsilon$,
the short geodesics $\delta$ and $\gamma$
correspond to disjoint components $X(\delta)$ and $X(\gamma)$ of the
thin part of $X$.
The total mass of $|\phi_\delta|$ in $X(\gamma)$ is $O(1)$.
Now $a_\gamma \phi_\gamma$ is chosen
to cancel the residue of $\phi_\delta$ in $X(\gamma)$, so we also
have $\|a_\gamma \phi_\gamma\|_T = O(1)$.
Since $\|\phi_\delta\|_T \bab \ell_\delta(X)^{-1}$, we obtain
the bound above.
\enddemo

\proclaim{Theorem}
\label{thm:Ceps}
If $\Res_\gamma(\phi)=0$ for all geodesics with
$\ell_\gamma(X)<\varepsilon${\rm ,} then we have
$$
        \|\phi\|_\WP \le C(\varepsilon) \,\|\phi\|_T .
$$
\endproclaim

\demo{Proof}
Let $X_r$, $r=\varepsilon/2$, denote the subset of $X$ with
hyperbolic injectivity radius less than $r$.
Since the area of $X$ is $2\pi|\chi(S)|$,
the thick part $X-X_r$ can be covered by
$N(r)$ balls of radius $r/2$.  The $L^1$-norm of
$\phi$ on a ball $B(x,r)$ controls its $L^2$-norm
on $B(x,r/2)$, so we have: 
$$
        \int_{X-X_r} \rho^{-2} |\phi|^2
        \EQ O(\|\phi\|_T^2) .
$$

It remains to control the $L^2$-norm of $\phi$ over the thin part $X_r$.  
For $\varepsilon$ sufficiently small,
every component of $X_r$ is either a horoball neighborhood
of a cusp or a collar neighborhood $X_r(\gamma)$ of a geodesic
$\gamma$ with $\ell_\gamma(X) < \varepsilon$.

To bound the integral of $\rho^{-2}|\phi|^2$ over 
a collar $X_r(\gamma)$, identify the covering space 
$X_\gamma \arrow X$ with $A(R)$
as before, and note that (for small $\varepsilon$) we have
$X_r(\gamma) \subset \pi(A(T))$ with $T=\varepsilon_0 R$.
Since $\pi|A(T)$ is injective we have
$$
        \int_{X_r(\gamma)} \rho^{-2} |\phi|^2
        \le \int_{A(T)} \rho^{-2} |\pi^*\phi |^2.
$$
Now because $\Res_\gamma(\phi)=0$, we can use the Laurent
expansion on $A(R)$ to write
$$
        \pi^*\phi 
        \EQ 
        z f(z) \frac{dz^2}{z^2} + 
        \frac{1}{z} \,g \left( \frac{1}{z} \right) \frac{dz^2}{z^2} 
        \EQ
        F+G,
$$
where $f(z)$ and $g(z)$ are holomorphic on $\Delta(R)=\{z \st |z|<R\}$.
Then $F$ and $G$ are orthogonal in $L^2(A(R))$, so we have
$$
        \int_{A(T)} \rho^{-2} |\pi^*\phi |^2
        = \int_{A(T)} \rho^{-2} (|F|^2 + |G|^2) .
$$
The inclusion $A(R) \subset \Delta^*(R)$ contracts
the hyperbolic metric, so to obtain an upper bound
on the integral above
we can replace $\rho(z)$ 
with $\rho_{\Delta(R)^*}(z) = 1/|z \log (R/|z|)|$.
Moreover $|f(z)|^2$ is subharmonic, so its mean over the
circle of radius $t$ is an increasing function of $t$.
Combining these facts, we see that
\begin{eqnarray*}
        \int_{A(T)} \rho^{-2} |F|^2 
        & = &
        \int_{A(T)} \frac{|zf(z)|^2}{|z|^4\rho^2(z)} 
        \,|dz|^2 
         \le  
        \int_{\Delta(T)}
        |f(z)|^2 \, |\log (R/|z|)|^2 \, |dz|^2 \\
        & = &
        \int_0^T 
        t (\log (R/t))^2 
        \int_0^{2\pi}
        |f(te^{i \theta})|^2 \,d\theta \, dt \\
        & \le &   
        2 \pi T^2 |\log \varepsilon_0|^2  \sup_{S^1(T)} |f(z)|^2 
         =  O \left( \sup_{S^1(T)} |zf(z)|^2 \right).
\end{eqnarray*}
Applying a similar argument to $|G|^2$, we obtain
$$
        \int_{A(T)} \rho^{-2} |\pi^*\phi|^2
        = O \left(\sup_{S^1(T)} |z g(z)|^2+|z f(z)|^2 \right).
$$

Without loss of generality we may assume 
$\sup_{S^1(T)} |f(z)| \ge \sup_{S^1(T)} |g(z)|$.
Since $\rho \bab |dz|/|z|$ on $S^1(T)$, we then have
$$
        \sup_{S^1(T)} \frac{|\pi^*\phi|}{\rho^2}
        \bab \sup_{S^1(T)} |z f(z) + g(1/z)/z|
        \bab \sup_{S^1(T)} |z f(z)| .
$$
Now $\pi(S^1(T))$ is contained in the thick part $X-X_r$, so we may 
conclude that
$$
        \int_{X_r(\gamma)} \rho^{-2} |\phi|^2
        = O \left( \sup_{X-X_r} \rho^{-4} |\phi|^2 \right).
$$
But the sup-norm of $\phi$ in the thick part is controlled by
its $L^1$-norm, so finally we obtain
$$
        \int_{X_r(\gamma)} \rho^{-2} |\phi|^2
        = O(\|\phi\|_T^2 ) .
$$

The bound on the $L^2$-norm of $\phi$ over the cuspidal components
of the thin part $X_r$ is similar,
using the fact that $\phi$ has at worst simple poles at the cusps.
Since the number of components of $X_r$ is bounded in terms
of $|\chi(S)|$, we obtain
$\int_{X_r} \rho^{-2} |\phi|^2 = O(\|\phi\|_T^2)$, completing
the proof.
\enddemo

{\it Remark}.
Masur has shown the Weil-Petersson metric extends
to $\closure{\cM}_{g,n}$, using a construction similar
to the thick-thin decomposition to
trivialize the cotangent bundle 
of $\cM_{g,n}$ near a curve with nodes
\cite{Mas}.  

\section{The $1/\ell$ metric}
\label{sec:comparison}

In this section we turn to
the K\"ahler metric $\gol$ on Teichm\"uller space, and
show it is comparable to the Teichm\"uller metric.

Recall that a positive $(1,1)$-form $\omega$ 
on $\Teich(S)$ determines
a Hermitian metric $g(v,w)=\omega(v,iw)$,
and $g$ is K\"ahler if $\omega$ is closed.
We say $g$ is {\it comparable} to the Teichm\"uller metric
if we have $\|v\|_T^2 \bab g(v,v)$ for all $v$ in the tangent
space to $\Teich(S)$.

\proclaimtitle{K\"ahler $\bab$ Teichm\"uller}
\proclaim{Theorem}
\label{thm:comparable}
Let $S$ be a hyperbolic surface of finite volume{\rm .}
Then for all $\varepsilon>0$ sufficiently small{\rm ,}
there is a $\delta>0$ such that the $(1,1)$\/{\rm -}\/form
\begin{equation}
\label{eq:wol}
                \wol =
                \omega_\WP - i\delta 
                \sumg
                \del \delbar \Log \frac{\varepsilon}{\ell_\gamma}
\end{equation}
defines a K{\rm \"{\it a}}hler metric $\gol$ on $\Teich(S)$ that is comparable to 
the Teichm{\rm \"{\it u}}ller metric{\rm .}
\endproclaim

Since the Teichm\"uller metric is complete we have:

\proclaimtitle{Completeness}

\proclaim{{C}orollary}
\label{cor:complete}
The metric $\gol$ is complete{\rm .}
\endproclaim 

\demo{Notation}
To present the proof of Theorem \ref{thm:comparable},
let $N = 3|\chi(S)|/2 + 1$
be a bound on the number of terms in the expression for $\wol$,
and let
$$
        \psi_\gamma =  \del \log \ell_\gamma =
                \frac{\del \ell_\gamma}{\ell_\gamma}; 
$$
 we then have
\begin{equation}
\label{eq:psi}
        |\psi_\gamma(v)|^2 = 
        \frac{i}{2} 
        \frac{\del\ell_\gamma \wedge \delbar \ell_\gamma}{\ell_\gamma^2}
        (v,iv) .
\end{equation}
\enddemo

\proclaim{Lemma} 
\label{lem:hermit}
There is a Hermitian metric $g$ of the form
\begin{equation}
\label{eq:gvv}
        g(v,v) =
        A(\varepsilon) \: \|v\|_\WP^2 \:+\:
        B \sumg |\psi_\gamma(v)|^2  
\end{equation}
such that $\|v\|_T^2 \le g(v,v) \le O(\|v\|_T^2)$
for all $\varepsilon>0$ sufficiently small{\rm .}
\endproclaim 

\demo{Proof}
By Propositions \ref{prop:WP:T} and \ref{thm:ell},
we have $\|v\|_\WP = O(\|v\|_T)$ and 
$\|\psi_\gamma(v)\| \le 2 \|v\|_T$,
and there are at most $N$ terms in the sum (\ref{eq:gvv}),
so $g(v,v) \le O(\|v\|_T^2)$.

To make the reverse comparison for a given
$v \mem \T_X\Teich(S)$, pick $\phi \mem Q(X)$ with
$\|\phi\|_T=1$ and $\phi(v)=\|v\|_T$.  
So long as $\varepsilon > 0$ is sufficiently small, 
we can apply the thick-thin decomposition for quadratic differentials
(Theorem \ref{thm:thick:thin})
to obtain
\begin{equation}
\label{eq:phisum}
        \phi = \phi_0 + \sumg
                a_\gamma \psi_\gamma  
\end{equation}
with $\Res_\gamma(\phi_0)=0$ and with $\|\psi_\gamma\|_T \ge 1$.
(Recall from Theorem \ref{thm:ell} that $\phi_\gamma$ and $\psi_\gamma$
are proportional, and that $\|\psi_\gamma\|_T \arrow 2$ as
$\varepsilon \arrow 0$.)

By Theorem \ref{thm:thick:thin}
each term on the right in (\ref{eq:phisum}) has Teichm\"uller norm 
$O(\|\phi\|_T) = O(1)$.
Since the residues of $\phi_0$ along the short geodesics 
vanish, the Teichm\"uller and Weil-Petersson norms
of $\phi_0$
are comparable, with a bound depending on $\varepsilon$ (Theorem \ref{thm:Ceps}).
Therefore we have $|\phi_0(v)| \le D(\varepsilon) \|v\|_\WP$.
Since we have $\|\psi_\gamma\|_T \ge 1$ and 
$\|a_\gamma \psi_\gamma\|_T = O(1)$,
we also have $|a_\gamma| \le E$,
where $E$ is independent of $\varepsilon$.
So from (\ref{eq:phisum}) we obtain
$$
        \phi(v) =
        \|v\|_T = \phi(v) 
        \le D(\varepsilon) \: \|v\|_\WP +
                E \sumg |\psi_\gamma(v)| .
$$
There are at most $N$ terms in the sum above, so we have
$$
        \|v\|^2_T \LE
        N D(\varepsilon)^2 \|v\|^2_\WP +
                N E^2 \sumg |\psi_\gamma(v)|^2 .
$$
Setting $A(\varepsilon)=N D(\varepsilon)^2$ and $B=NE^2$, 
from (\ref{eq:gvv}) we obtain $\|v\|_T^2 \le g(v,v)$.
\enddemo

\proclaim{{C}orollary}
\label{cor:hermit}
For $\varepsilon > 0$ sufficiently small{\rm ,} we have
$$
        \|v\|_T^2 \bab
        \: \|v\|_\WP^2 \:+\:
        \sumg |(\del \log \ell_\gamma)(v)|^2  
$$
for all vectors $v$ in the tangent space to $\Teich(S)${\rm .}
\endproclaim 

Next we control the terms in (\ref{eq:wol}) coming
from geodesics of length near~$\varepsilon$.

\proclaim{Lemma} 
\label{lem:medium}
For $\varepsilon < \ell_\delta(X) < 2\varepsilon$ we have
$$
        |\psi_\delta(v)|^2
        \le D(\varepsilon) \|v\|_\WP^2 +
        O \left( \varepsilon \sumg |\psi_\gamma(v)|^2 \right)
$$
for any tangent vector $v \mem \T_X \Teich(S)${\rm .}
\endproclaim 

\demo{Proof}
By Theorem \ref{thm:gamma:medium} we have
$\psi_\delta = \psi_0 + \sumg a_\gamma \psi_\gamma$
with $a_\gamma = O(\ell_\delta(X)) = O(\varepsilon)$, and with
$\|\psi_0\|_T \le C(\varepsilon) \|\psi_\delta\|_\WP$
by Theorem \ref{thm:Ceps}.  Evaluating this sum
on $v$, we obtain the lemma.
\enddemo

{\it Proof of Theorem} \ref{thm:comparable} 
        (K\"ahler $\bab$ Teichm\"uller). 
Consider the $(1,1)$-form
\begin{equation}
\label{eq:w}
        \omega = 
        (F(\varepsilon) + A(\varepsilon)) \omega_\WP
        - B \sumd \frac{i}{2} \:
                \del \delbar \Log \frac{2\varepsilon}{\ell_\delta} ,
\end{equation}
where 
$$
        F(\varepsilon)=16NBD(\varepsilon) \,\sup_{[1,2]} |\Log''(x)|,
$$
and where 
$A(\varepsilon)$, $B$ and $D(\varepsilon)$ come from the lemmas above.
Then $\omega$ and $\wol$ are of the same form
(up to scaling and replacing $\varepsilon$ with $\varepsilon/2$),
so to prove the Theorem it suffices to show 
$g(v,v) = \omega(v,iv)$ is comparable
to the Teichm\"uller metric.

Let $v \mem \T_X \Teich(S)$ be a vector with $\|v\|_T=1$.
To begin the evaluation of $g(v,v)$, we compute
\begin{equation}
\label{eq:dd}
        \del \delbar \Log \frac{2\varepsilon}{\ell_\delta} 
        = 2 \varepsilon 
        \left (\Log' \frac{2\varepsilon}{\ell_\delta} \right)
        \del \delbar \left( \frac{1}{\ell_\delta} \right)
        \;+\;
        \frac{4 \varepsilon^2}{\ell_\delta^2}
        \left(
                \Log'' \frac{2\varepsilon}{\ell_\delta} 
        \right)
        \frac{\del\ell_\delta \wedge \delbar \ell_\delta}{\ell_\delta^2} 
        \cdot \hskip.35in
\end{equation}
By Theorem \ref{thm:ah}, the function $1/\ell_\delta$ is
almost pluriharmonic; more precisely,
$$
        \left|\del \delbar 
                \left( \frac{1}{\ell_\delta} \right)
        (v,iv) \right| = O(1).
$$
Since $\Log'(x)$ is bounded, the term in (\ref{eq:dd}) involving
$\del\delbar (1/\ell_\delta)$ is $O(\varepsilon)$.
Using (\ref{eq:psi}) we then obtain
\begin{equation}
\label{eq:diff}
        \frac{i}{2}
        \left( \del \delbar \Log \frac{2\varepsilon}{\ell_\delta} \right)
        (v,iv)
        \EQ
        \frac{4 \varepsilon^2}{\ell_\delta^2}
        \left( \Log'' \frac{2\varepsilon}{\ell_\delta} \right)
        |\psi_\delta(v)|^2  
        + O(\varepsilon)  .
\end{equation}

Using expression (\ref{eq:w}) to compute $\omega(v,iv)$,
we obtain a sum of terms like that above, with $\ell_\delta(X) < 2\varepsilon$.
If $\ell_\delta(X) < \varepsilon$, then
$\Log''(2\varepsilon/\ell_\gamma)=
        \log''(2\varepsilon/\ell_\gamma)
        =-\ell_\gamma^2/(4\varepsilon^2)$;
  hence
$$
        -\frac{i}{2}
        \left( \del \delbar \Log \frac{2\varepsilon}{\ell_\delta} \right)
        (v,iv)
        = |\psi_\delta(v)|^2 + O(\varepsilon) .
$$
On the other hand, if $\varepsilon \le \ell_\delta(X) < 2\varepsilon$,
then from (\ref{eq:diff}) and Lemma \ref{lem:medium} we obtain:
\begin{eqnarray*}
        \left|
        \frac{i}{2}
        \left( \del \delbar \Log \frac{2\varepsilon}{\ell_\delta} \right)
        (v,iv)
        \right|
        &  \le &
        16 |\psi_\delta(v)|^2  \;\sup_{[1,2]} |\Log''(x)| 
        + O(\varepsilon)\\
        & \le & 
        \frac{F(\varepsilon)}{NB} \|v\|_\WP^2
        + O \left( \varepsilon \sumg |\psi_\gamma(v)|^2 \right)
        + O(\varepsilon)  .
\end{eqnarray*}
Applying these two bounds to $g(v,v)=\omega(v,iv)$, we obtain:
\begin{eqnarray*}
        g(v,v) & \ge &
        A(\varepsilon) \|v\|^2_\WP  +
        B \sumg |\psi_\gamma(v)|^2 \\
        && +\;
        B \sumdg \left( 
                \frac{F(\varepsilon)}{NB} \|v\|^2_\WP -
        \left|
        \frac{i}{2}
        \left( \del \delbar \Log \frac{2\varepsilon}{\ell_\delta} \right)
        (v,iv)
        \right|
        \right)
        + O(\varepsilon) 
        \\
        & \ge & 
        A(\varepsilon) \|v\|^2_\WP  +
        B (1+O(\varepsilon)) \sumg |\psi_\gamma(v)|^2 +O(\varepsilon) .
\end{eqnarray*}
By Lemma \ref{lem:hermit} we then have:
$$
        g(v,v) \ge \|v\|^2_T + O(\varepsilon) = 1+O(\varepsilon) .
$$
Thus $\|v\|_T^2 = O(g(v,v))$ when $\varepsilon$ is small enough.
The reverse comparison, $g(v,v) = O(\|v\|_T^2)$, follows the same lines
as Lemma \ref{lem:hermit}.
\hfill\qed

\vglue-3pt
\section{Quasifuchsian reciprocity}
\label{sec:reciprocity}
\vglue-3pt

Let $S$ be a hyperbolic surface of finite area
with quasifuchsian space $QF(S) = \Teich(S) \times \Teich(\Sbar)$.
In this section we define a map
\smallbreak
\centerline{$q : \T QF(S) \arrow \cx,$}
\smallbreak\noindent 
providing a natural bilinear pairing $q(\mu,\nu)$ on
$M(X) \times M(Y)$ for each $(X,Y) \mem QF(S)$.
The symmetry of this pairing (explained below)
will play a key role
in the discussion of a bounded primitive for the
Weil-Petersson metric in the next section.

To define $q$, recall that
a small change in the conformal structure on
$X$ determines a change in the projective structure on $Y$,
by the derivative of the Bers embedding
$$
        D\beta_Y : \T_X \Teich(S) \arrow P(Y) .
$$
Since $S$ has finite area, we also have
$$
        P(Y) \isom Q(Y) \isom \T^*_Y \Teich(\Sbar) .
$$
Thus for $(\mu,\nu) \mem M(X) \times M(Y)$, we can define
the {\it quasifuchsian pairing} 
$$
        q(\mu,\nu) = \brackets{D\beta_Y(\mu),\nu}
$$
by evaluating the cotangent vector $D\beta_Y(\mu)$ on the tangent
vector $[\nu] \mem {\rm T}_Y\Teich(S)$.  This pairing only depends on
the equivalence class represented by 
$(\mu,\nu)$ in $\T_{(X,Y)}QF(S)$.

Interchanging the roles of $X$ and $Y$, we obtain
a similar pairing from the Bers embedding
$$
        \beta_X : \Teich(\Sbar) \arrow P(X) .
$$
The main result of this section is that these pairings
are equal.

\proclaimtitle{Quasifuchsian reciprocity}
\proclaim{Theorem}
\label{thm:reciprocity}
For any $(\mu,\nu) \mem \T_{(X,Y)} QF(S)${\rm ,} we have
$$
        q(\mu,\nu) \EQ
        \int_Y (D\beta_Y(\mu)) \cdot \nu 
        \EQ
        \int_X (D\beta_X(\nu)) \cdot \mu . 
$$
\endproclaim

\noindent This theorem says $q$ is {\it symmetric}, meaning
$q \compos D\rho = q$,
where 
$$
        D\rho : \T QF(S) \arrow \T QF(S)
$$
is the derivative of
the involution $\rho(X,Y)=(\Ybar,\Xbar)$.

\demo{Proof}
Recall that $1/(\pi z)$ is the fundamental solution to
the $\delbar$ equation on $\chat$.
Thus a solution to the infinitesimal Beltrami equation 
$\delbar v = \mu$ is given by 
$$
        v(z) \frac{\del}{\del z} \EQ
        \left(
        \frac{1}{\pi} \int_\chat \frac{\mu(w)}{(z-w)} \, |dw|^2 
        \right)
        \frac{\del}{\del z} \cdot
$$
Since $1/(\pi z)''' = -6/(\pi z^4)$,
outside the support of $\mu$ the vector field
$v$ is holomorphic with
infinitesimal Schwarzian derivative given by
$$
        \phi = v'''(z) \,dz^2 = K*\mu
$$
where $K$ is the kernel
$$
        K = -\frac{6}{\pi (z-w)^4}
        \,dz^2\,dw^2.
$$
This kernel on $\chat \times \chat$ is natural 
and symmetric, in the sense that:
\begin{itemize}
        \item 
        $(\gamma \times \gamma)^*K = K$
for any M\"obius transformation $\gamma \mem \Aut(\chat)$; and
        \item
        $\iota^*K = K$ where $\iota(w,z)=(z,w)$.
\end{itemize}
The symmetry of $q$ will come from the symmetry of $K$.

To compute the derivative of Bers' embedding $\beta_Y$,
let us regard $\mu \mem M(X)$ and $\phi = D\beta_Y(\mu) \mem Q(Y)$
as $\Gamma(X,Y)$-invariant forms on $\Omega(X,Y)$.
Then we have
\begin{eqnarray}
\noalign{\vskip4pt}
\label{eq:Bers}
        &&D\beta_Y(\mu) =
        \phi(z) \, dz^2 =
        \left(
        -\frac{6}{\pi} \int_\chat \frac{\mu(w)}{(z-w)^4} \, |dw|^2 
        \right) \,dz^2 . \\
\noalign{\vskip-6pt}
&&  \nonumber
\end{eqnarray}
(See \cite{Ber}, 
\cite[\S 5.7]{Gd}.)

Now consider the kernel on $\chat \times \chat$
given by
\begin{eqnarray*}
\noalign{\vskip4pt}
&&
        K_0 = \sum_{\gamma \mem \Gamma(X,Y)} (\gamma,\id)^*K .
\\
\noalign{\vskip-6pt}
&&  
\end{eqnarray*}
Since $K$ was already invariant under the diagonal 
action of $\Aut(\chat)$, the kernel $K_0$ is invariant under
$\Gamma(X,Y) \times \Gamma(X,Y)$, and so it descends to
a form on $X \times Y$.  We then have 
\begin{eqnarray*}
\noalign{\vskip4pt}
&&
        q_0(\mu,\nu)
        = \brackets{D\beta_Y(\mu),\nu}
        = \int_{X \times Y} K_0(w,z) \mu(w) \nu(z) 
                \,|dw|^2 \,|dz|^2 .
\\
\noalign{\vskip-6pt}
&&  
\end{eqnarray*}
The reverse pairing is given similarly by
\begin{eqnarray*}
\noalign{\vskip4pt}
&&
        q_1(\mu,\nu)
        = \brackets{D\beta_X(\nu),\mu}
        = \int_{Y \times X} K_1(w,z) \nu(w) \mu(z) 
                \,|dw|^2 \,|dz|^2 ,
\\
\noalign{\vskip-6pt}
&&  
\end{eqnarray*}
where the form $K_1$ on $Y \times X$ is given upstairs by
$$
        K_1 = \sum_{\gamma \mem \Gamma(X,Y)} (\id,\gamma)^*K .
$$
By symmetry of $K$, $K_0 = \pi^*(K_1)$ for the
natural map $\pi : X \times Y \arrow Y \times X$;
therefore $q_0(\mu,\nu)=q_1(\mu,\nu)$, establishing reciprocity.
\enddemo

\section{The Weil-Petersson form is $d$(bounded)}
\label{sec:dbounded}

In this section we construct an explicit $1$-form
$\theta_\WP$ on $\Teich(S)$ such that $d(i\theta_\WP)=\omega_\WP$.
Using this $1$-form and Nehari's bound for
univalent functions, we then show the 
K\"ahler metrics corresponding to
$\omega_\WP$ and $\wol$ are both $d$(bounded).

Recall $\sigqf(X,Y) \mem \Proj_X(S)$ is the projective structure
on $X$ coming from the quasifuchsian uniformization of $X \cup Y$.

\proclaimtitle{Quasifuchsian primitive}
\proclaim{Theorem}
\label{thm:lag}
Fix $Y \mem \Teich(\Sbar)${\rm , }
and let $\theta_\WP$ be the $(1,0)$\/{\rm -}\/form on $\Teich(S)$ given by
\begin{eqnarray*}
        \theta_\WP(X) & = & \sigf(X) - \sigqf(X,Y) \\
        & = & -\beta_X(Y) .
\end{eqnarray*}
Then $i\theta_\WP$ is a primitive for the Weil\/{\rm -}\/Petersson
K{\rm \"{\it a}}hler form{\rm ;} that is{\rm ,} $d(i\theta_\WP)=\omega_\WP${\rm .}
\endproclaim

A key ingredient in the proof is:

\proclaim{Theorem}
For any $Y_0, Y_1 \mem \Teich(\Sbar)${\rm ,} we have
$$
        d(\sigqf(X,Y_1)-\sigqf(X,Y_0)) = 0 
$$
as a $2$\/{\rm -}\/form on $\Teich(S)$\/{\rm .}\/
\endproclaim

\demo{Proof}
Let $Y_t$ be a smooth path in $\Teich(\Sbar)$ joining
$Y_0$ to $Y_1$, and let 
$$
        \theta_t(X)= \sigqf(X,Y_t)-\sigqf(X,Y_0) .
$$
We will show $\theta_1=dF$ for an explicit function
$F : \Teich(S) \arrow \cx$.

Let $\nu_t = dY_t/dt \mem \T_{Y_t} \Teich(\Sbar)$.
The quasifuchsian uniformization of $(X,Y_t)$ determines
a Schwarzian quadratic differential $\beta_{Y_t}(X)$ on $Y_t$.
Let
$$
        f_t(X) = \brackets{\sigqf(X,Y_t)-\sigf(Y_t),\nu_t}
                = \brackets{\beta_{Y_t}(X),\nu_t} .
$$

We claim $d f_t = \partial \theta_t/\partial t$ as $1$-forms on $\Teich(S)$.
Indeed, by quasifuchsian reciprocity, 
for any $\mu \mem M(X)$ we have 
$$
        df_t(\mu) =
        \brackets{D\beta_{Y_t}(\mu),\nu_t} =
        \brackets{D\beta_{X}(\nu_t),\mu} =
        \left\langle
        \frac{\partial \sigqf(X,Y_t)}{\partial t},\mu \right\rangle
        = \frac{\partial \theta_t}{\partial t}(\mu) \cdot
$$
Set $F(X) = \int_0^1 f_t(X) \,dt$.
Then, since $\theta_0=0$,  
$$
        \theta_1 = 
        \int_0^1 \frac{\partial \theta_t}{\partial t} \,dt
        =
        d \int_0^1 f_t \,dt
        =dF,
$$
and thus $d\theta_1=d^2F = 0$.
\enddemo

{\it Proof of Theorem} \ref{thm:lag} (Quasifuchsian primitive). 
Let us compute the $2$-form $d\theta_\WP$ at $X_0 \mem \Teich(S)$.
By the preceding result, we may freely modify the
choice of $Y$ without changing $d\theta_\WP$.
Setting $Y=\Xbar_0$ we obtain
\begin{eqnarray*}
\noalign{\vskip4pt}
&&
        \theta_\WP =
        \sigqf(X,\Xbar) - \sigqf(X,\Xbar_0) .
\\
\noalign{\vskip-6pt} 
&&
\end{eqnarray*}
Now $\sigqf(X,\Ybar)$ is holomorphic in $X$ and antiholomorphic
in $Y$, so to compute the $(2,0)$ part of $d\theta_\WP$ we can
replace $\Xbar$ by $\Xbar_0$ in $\sigqf(X,\Xbar)$; then
we obtain:
\begin{eqnarray*}
\noalign{\vskip4pt}
&&
        \del\theta_\WP = \del(\sigqf(X,\Xbar_0) - \sigqf(X,\Xbar_0)) = 0 .
\\
\noalign{\vskip-6pt} 
&&
\end{eqnarray*}
As for the $(1,1)$-part, we can similarly replace $X$ by $X_0$
in $-\sigqf(X,\Xbar_0)$ to obtain:
$$
        \delbar\theta_\WP =
        \delbar(\sigqf(X_0,\Xbar) - \sigqf(X_0,\Xbar_0)) 
        = \delbar \beta_{X_0}(\Xbar) .
$$
To complete the proof, it suffices to show
$$
        (\delbar\theta_\WP)(\mu,i\mu) = -i \|\mu\|_\WP^2
$$
for all $[\mu] \mem \T_{X_0} \Teich(S)$.

Since $S$ has finite hyperbolic area, any tangent direction
to $\Teich(S)$ at $X_0$ 
is represented by a harmonic Beltrami differential
\begin{eqnarray*}
\noalign{\vskip4pt}
&&
        \mu = \rho^{-2} \phibar, 
\\
\noalign{\vskip-6pt} 
&&
\end{eqnarray*}
where $\phi \mem Q(X_0)$ and $\|\mu\|_\WP = \|\phi\|_\WP$.
Let $\mubar \mem M(\Xbar_0)$ be the corresponding
conjugate vector tangent to $\Teich(\Sbar)$ at $\Xbar_0$.
Since $\theta_\WP(X_0) = 0$, we have
\begin{eqnarray*}
        (\delbar\theta_\WP)(\mu,i\mu)  & = & 
        \brackets{D\beta_{X_0}(\mubar),i\mu} -
        \brackets{D\beta_{X_0}(-i\mubar),\mu} \\ \noalign{\vskip4pt}
        & = &
        2 i \,\brackets{D\beta_{X_0}(\mubar),\mu} .
\end{eqnarray*}
The evaluation of $D\beta_{X_0}(\mubar)$ is a standard
calculation of the derivative of the Bers embedding at
the origin.  Namely writing $X_0=\half/\Gamma_0$, we can
interpret $\phi$ and $\mubar$ as $\Gamma_0$-invariant
forms on $\half$ and $\lowerhalf = -\half$, respectively.
Then $\mubar(w) = \phi(\wbar)$, and by (\ref{eq:Bers}),
$$
        D\beta_{X_0}(\mubar) = \psi(z) \, dz^2 =
        \left(
        -\frac{6}{\pi} \int_\lowerhalf
        \frac{\rho^2(\wbar) \phi(\wbar)}{(z-w)^4} \, |dw|^2 
        \right) \,dz^2 .
$$
A well-known reproducing formula
(see \cite[\S 5.7]{Gd},
\cite[(5.2)]{Ber}) gives
$\psi = (-1/2) \phi$.  Therefore we have
$$
        (\delbar\theta_\WP)(\mu,i\mu) = 
        2i \,\brackets{(-1/2)\phi,\mu} 
        = -i \int_X \phi \rho^{-2} \phibar
        = -i \,\|\mu\|_\WP^2 ,
$$
which completes the proof.
\hfill\qed

\proclaimtitle{$d$(bounded)}
\proclaim{{C}orollary} 
\label{cor:dbounded}
The K{\rm \"{\it a}}hler form of the $1/\ell$ metric
on Teichm{\rm \"{\it u}}ller space is $d${\rm (}\/bounded\/{\rm ), } with primitive
\begin{equation}
\label{eq:tol}
        \tol \EQ
        \theta_\WP - \delta 
                \sumg
                \del \Log \frac{\varepsilon}{\ell_\gamma} 
\end{equation}
satisfying $d(i\tol) = \wol${\rm .}
\endproclaim 

\demo{Proof}
By (\ref{eq:wol}) we have $\wol = d(i\tol)$.
To bound the first term, we note that
$-\theta_\WP(X) = \sigqf(X,Y) - \sigf(X)$ is the 
Schwarzian derivative of the univalent map 
$$
        f_{X,Y} : \half \arrow  \widetilde{X} \subset \Omega(X,Y) \subset \chat,
$$
so by Nehari's bound (Theorem~\ref{thm:Nehari}) we have
$\|\theta_\WP(X)\|_\infty < 3/2$.  Therefore
\begin{equation}
\label{eq:L1}
        \|\theta_\WP\|_T 
        = \int_X \frac{|\phi|}{\rho^2} \rho^2
        \le 3\pi |\chi(S)|,
\end{equation}
since $\int \rho^2 = 2\pi |\chi(S)|$ by
Gauss-Bonnet.

For the remaining terms in (\ref{eq:tol}) we appeal to
Theorem \ref{thm:ell}, which gives the
bound $\|\del \ell_\gamma\|_T \le 2 \ell_\gamma$.
The latter implies:
$$
        \left\|
        \del \Log \frac{\varepsilon}{\ell_\gamma} \right\|_T
        \EQ \left| \frac{\varepsilon}{\ell_\gamma^2}
        \Log' \frac{\varepsilon}{\ell_\gamma} \right|
        \cdot \|\del \ell_\gamma\|_T
        \LE
        \left| \frac{\varepsilon}{\ell_\gamma}
        \Log' \frac{\varepsilon}{\ell_\gamma} \right|
        \EQ O(1) .
$$
Putting these bounds together, we have
$\|\tol\|_T = O(1)$.
Since the $1/\ell$ metric is comparable to the
Teichm\"uller metric (Theorem~\ref{thm:comparable}),
we have shown that $\wol$ is $d$(bounded).
\enddemo  

The $L^2$-version of (\ref{eq:L1}) gives
$\|\theta_\WP\|_\WP^2 \le (9\pi/2) |\chi(S)|$, so we also have:

\proclaim{{C}orollary}
The K{\rm \"{\it a}}hler form $\omega_\WP$ is 
$d${\rm (}\/bounded\/{\rm )} for the Weil\/{\rm -}\/Petersson metric{\rm .}
\endproclaim 

\section{Volume and curvature of moduli space}
\label{sec:volume}

In this section we prove that $(\cM(S),\gol)$ has bounded
geometry and finite volume,
completing the proof of the K\"ahler hyperbolicity of
moduli space.

\proclaimtitle{Finite volume}
\proclaim{Theorem} 
The metric $\gol$ descends to the moduli space 
$\cM(S)${\rm ,} and $(\cM(S),\gol)$ has finite volume{\rm .}
\endproclaim

\demo{Proof}
By its definition (\ref{eq:wol}), the metric $\gol$ is
invariant under the action of the mapping class group
$\Mod(S)$ on $\cM(S)$, so it descends to a metric on moduli space.

Let $n=\dim_{\cx} \cM(S)$, and let
$\closure{\cM(S)}$ be the Deligne-Mumford compactification 
of $\cM(S)$.
Consider a stable curve
$Z \mem \closure{\cM(S)}$  with $k$ nodes.
Then $Z$ has a neighborhood
$U$ in $\closure{\cM(S)}$ satisfying
\begin{eqnarray*}
        (U,Z) & \isom & (\Delta^n,0)/G \;\;\;\;{\rm and}\\
        U \cap \cM(S) & \isom & ((\Delta^*)^k \times \Delta^{n-k})/G,
\end{eqnarray*}
where $G$ is a finite group.
(Compare \cite[\S 3]{Wol3}.)

A small neighborhood of $(0,0)$ has finite volume in the
Kobayashi metric on $(\Delta^*)^k \times \Delta^{n-k}$,
and inclusions contract the Kobayashi metric,
so there is a neighborhood $V$ of $Z$ in $\closure{\cM(S)}$ with
$\vol(V \cap \cM(S)) < \infty$ in the 
Kobayashi $=$ Teichm\"uller  metric on $\cM(S)$.  
By compactness of $\closure{\cM(S)}$, the Teichm\"uller volume
of $\cM(S)$ is finite.

Since the Teichm\"uller metric
is comparable to $\gol$, $\cM(S)$ also has finite volume
in the $1/\ell$ metric.
\enddemo

\proclaimtitle{Bounded geometry}
\proclaim{Theorem}
The sectional curvatures of the metric $\gol$ 
are bounded above and below over $\Teich(S)${\rm ,} and the
injectivity radius of $\gol$ is uniformly bounded below{\rm .}
\endproclaim

\demo{Proof}
We use the method of the proof of Theorem \ref{thm:ah}:
namely we realize $\Teich(S)$ as the locus $(X,\Xbar)$ in
$QF(S)$, and extend the
functions $\sigqf$ and $\ell_\gamma$ used in the
definition of $\gol$ to holomorphic functions on $QF(S)$.
Uniform bounds on these extensions then give
$C^\infty$ bounds on $\gol$.

Pick $X_0 \mem \Teich(S)$, and let $n=\dim_{\cx} \Teich(S)$.
By an elementary result in convex geometry, there is
an isomorphism $A : \cx^n \arrow P(\Xbar_0)$
such that $\|Az\|_\infty \bab |z|$ with constants
depending only on $n$.
Using the Bers embedding of $\Teich(S)$ into $P(\Xbar_0)$
and Theorem \ref{thm:Bers},
we can find an embedded polydisk
\smallbreak
\centerline{$
        f : (\Delta^n,0) \arrow (\Teich(S),X_0),
$}\smallbreak\noindent 
such that the Teichm\"uller
and Euclidean metrics are comparable on $\Delta^n$.
(Indeed we can take $f(z) = A(\alpha z)$ for
a suitable $\alpha > 0$.)

Let $X_s = f(s) \mem \Teich(S)$ and 
$Y_t = \Xbar_\tbar \mem \Teich(\Sbar)$;
then $(X_s,Y_t) \mem QF(S)$ is a holomorphic function of 
$(s,t) \mem \Delta^n \times \Delta^n$.
For any closed geodesic $\gamma$ on $S$, the complex length
satisfies $\Re\, \cL_\gamma(X_s,Y_t) > 0$, so
the holomorphic function
\smallbreak
\centerline{$
        L_\gamma(s,t) = \log \cL_\gamma(X_s,Y_t)
$}\smallbreak\noindent 
maps $\Delta^n \times \Delta^n$ into the strip
$\{z \st |\Im z| < \pi\}$.  By the Schwarz lemma,
all the derivatives of $L_\gamma(s,t)$ at $(0,0)$ of
order $\le k$ are bounded by a constant $C_k$
in the Euclidean metric.
Therefore the derivatives of 
\smallbreak
\centerline{$
        \log \ell_\gamma(X_s) = \log \cL_\gamma(X_s,\Xbar_s)
$}\smallbreak\noindent 
are also controlled at $s=0$. \pagebreak

Fixing $Z \mem \Teich(\Sbar)$, the holomorphic $1$-form
$$
        \tau(X,Y) = \sigqf(X,Y) - \sigqf(X,Z)
$$
is bounded in the Teichm\"uller metric on $QF(S)$
(by Nehari's bound, Theorem~\ref{thm:Nehari}). 
Hence $f^*\tau$ is bounded in the Euclidean metric,
and thus the derivatives of
$$
        (f^*\tau)(s) = 
        f^*(\sigf(X_s) - \sigqf(X_s,Z))
$$
at $s=0$ are also uniformly controlled.

Now by Corollary \ref{cor:dbounded}, if we set 
$$
        \theta \EQ
        f^*(\sigf(X_s) - \sigqf(X_s,Z))
        \;-\; \delta 
        \sumg
        \del \Log \frac{\varepsilon}{\ell_\gamma(X_s)} ,
$$
then the 1-form $i\theta$
is a primitive for the K\"ahler form of $g=f^*\gol$ on $\Delta^n$.
Since the derivatives of each term are controlled, 
and $g$ is comparable to the Euclidean metric, 
we find that the sectional curvatures of $g$
are bounded above and below at $s=0$.
(Indeed $g$ ranges in a 
compact family of metrics in the $C^\infty$ topology on $\Delta^n$.)
Thus the curvatures of $\gol$ are $O(1)$ at $X_0$, and hence
uniformly bounded over $\Teich(S)$.  

We then have curvature bounds on $g$ throughout $\Delta^n$,
from which it follows that the injectivity radius of $g$
at $0$ is bounded below.  Since $f$ is one-to-one,
the injectivity radius of $\gol$ on $\Teich(S)$ is also
bounded below.
\enddemo

{\it Proof of Theorem} \ref{thm:kh} (K\"ahler hyperbolic). 
The metric $h=\gol$ is comparable to the Teichm\"uller metric by
Theorem \ref{thm:comparable}, and therefore $\gol$ is
complete (Corollary \ref{cor:complete}).
Its symplectic form $\wol$ is $d$(bounded) by 
Corollary \ref{cor:dbounded}.
In this section we have shown that $(\cM(S),\gol)$ has finite
volume and bounded geometry, and therefore moduli space
is K\"ahler hyperbolic in the $\gol$ metric.
\phantom{out of time}\hfill\qed

\section{Appendix.  Reciprocity for Kleinian groups}

This appendix formulates a version of
quasifuchsian reciprocity for general Kleinian groups.
As an application, we sketch a proof of the 
Takhtajan-Zograf formula
$$
        d(\sigma_F(X) - \sigma_S(X)) = -i\omega_\WP .
$$

\demo{Kleinian groups}
Let $\Gamma \subset \Aut (\chat)$ be a finitely
generated Kleinian group, with domain of
discontinuity $\Omega$ and (possibly disconnected)
quotient Riemann surface $X=\Omega/\Gamma$.
By Ahlfors' finiteness theorem, $X$ has finite
hyperbolic area.

Let $\mu \mem M(X)$ be a Beltrami differential,
regarded as a $\Gamma$-invariant form on $\chat$,
vanishing outside $\Omega$.
Let $v$ be a quasiconformal vector field on the
sphere with $\delbar v = \mu$.  Then the
infinitesimal Schwarzian derivative
$$
        \phi_\mu = \frac{\del^3 v}{\del z^3} \,dz^2
$$
is a $\Gamma$-invariant quadratic differential,
holomorphic {\it outside} the support of $\mu$.

On the support of $\mu$, the third derivative of $v$
exists only as a distribution.  If, however, $\mu$
is sufficiently smooth --- for example, if
$\mu$ is a harmonic Beltrami differential
($\mu = \rho^{-2} \phibar$, $\phi \mem Q(X)$) ---
then $\phi_\mu$ is smooth on $\Omega$ and descends to
a quadratic differential
$$
        K(\mu) = 
        \phi_\mu \mem L^1(X,dz^2) .
$$
We refer to $\phi_\mu$ as the {\it projective distortion}
of $\mu$.
\enddemo

\proclaimtitle{Kleinian reciprocity}
\proclaim{Theorem}
Let $X=\Omega/\Gamma$ be the quotient Riemann surface
for a finitely generated Kleinian group $\Gamma${\rm ,} and
let $\mu, \nu \mem M(X)$ be a pair of 
sufficiently smooth Beltrami differentials{\rm .}
Then we have\/{\rm :}\/
$$
        \int_X \phi_\mu \, \nu = \int_X \phi_\nu \, \mu ,
$$
where $\phi_\mu, \phi_\nu \mem L^1(X,dz^2)$ give
the projective distortions of $\mu$ and $\nu${\rm .}
\endproclaim

Note that if $\Gamma$ is quasifuchsian,
with $X=X_1 \disjunion X_2$, then by taking
$\mu \mem M(X_1)$ and $\nu \mem M(X_2)$ we recover
Theorem \ref{thm:reciprocity} (quasifuchsian reciprocity).
The proof of the general version is essentially the same;
it turns on the symmetry of the kernel 
$K = -6 \,dz^2\,dw^2/(\pi (z-w)^4)$.

\demo{Schottky uniformization}
As an application, let $S$ be a closed surface of
genus $g \ge 2$, and fix a maximal collection
$(a_1,\ldots,a_g)$ of disjoint simple closed
curves on $S$, linearly independent in $H_1(S,\reals)$.
Then to each $X \mem \Teich(S)$ one can associate
a {\it Schottky group} $\Gamma_S$, such that
$X=\Omega_S/\Gamma_S$ and the curves $a_i$ lift
to $\Omega_S$.  Let $\sigma_S(X)$ be the 
projective structure on $X$ inherited from
$\Omega_S$.   \enddemo

\proclaimtitle{Takhtajan-Zograf}  
\proclaim{Theorem}
The difference between the
Fuchsian and Schottky projective structures gives
a $1$\/{\rm -}\/form on $\Teich(S)$ satisfying\/{\rm :}\/
$$
        d(\sigma_F(X) - \sigma_S(X)) = -i\omega_\WP .
$$
\endproclaim

\demo{Sketch of the proof}
Since $\sigma_S$ and $\sigqf$ are holomorphic,
we have
$$\delbar (\sigma_F-\sigma_S) = \delbar (\sigma_F-\sigqf)
        = -i\omega_\WP$$ (Theorem~\ref{thm:lag}).
Thus we just need to check that the  (2,0)-form
$\del (\sigma_F-\sigma_S)$ on $\Teich(S)$ vanishes.
To this end, let us represent the tangent space to
$\Teich(S)$ at $X$ by the space $\cH(X)$ of harmonic
Beltrami differentials.  Using the Fuchsian and Schottky
representations of $X$ as a quotient Riemann surface
for Kleinian groups $\Gamma_F$ and $\Gamma_S$, we obtain operators
$$
        K_F, K_S : \cH(X) \arrow L^1(X,dz^2)
$$
sending $\mu \mem \cH(X)$ to its projective distortion $\phi_\mu$.
(In the Fuchsian case, $\Omega_F/\Gamma_F = X \disjunion \Xbar$;
we set $\mu=0$ on $\Xbar$.)

We then compute:
$$
        \del (\sigma_F-\sigma_S)(\mu,\nu) \EQ
        \brackets{K_F(\mu)-K_S(\mu),\nu} -
        \brackets{K_F(\nu)-K_S(\nu),\mu},
$$
where $\brackets{\phi,\mu} = \int_X \phi \, \mu$.
The key to this computation is to observe that the
usual chain rule for the Schwarzian derivative,
$S(g \compos f) = Sf + f^*Sg$, holds even when
just one of $f$ and $g$ is conformal (assuming the
other is sufficiently smooth).

By reciprocity, the bracketed expressions above agree,
so $\del (\sigma_F-\sigma_S)=0$.
\phantom{walkingonsunshine}\hfill\qed\enddemo

\AuthorRefNames [Wol3]

\bye